\documentclass[conference]{IEEEtran}
\IEEEoverridecommandlockouts
\usepackage{lineno}
\usepackage[cmex10]{amsmath}
\usepackage{amscd,amsfonts,amsmath,amssymb,mathrsfs,pifont,stmaryrd,tipa}
\usepackage{array,cases,dsfont,graphicx,texdraw}
\usepackage{graphics}
\usepackage{epstopdf}
\usepackage{epsfig}
\usepackage{cite}
\usepackage{mathbbold}
\usepackage{stfloats}
\usepackage{subfig}
\usepackage{hyperref}
\usepackage{amsfonts}
\usepackage{amssymb}
\usepackage{graphicx}%
\usepackage{xcolor}
\newtheorem{theorem}{\bf{Theorem}}

\newtheorem{assumption}{Assumption}

\newtheorem{remark}{\bf{Remark}}

\allowdisplaybreaks
\begin{document}

\title{A Continuous-Time Nesterov Accelerated Gradient Method for Centralized and Distributed Online Convex Optimization}
\author{Chao Sun and Guoqiang Hu\thanks{This work was supported by Singapore Ministry of Education Academic Research Fund Tier 1 RG180/17(2017-T1-002-158). C. Sun and G. Hu are with the School of
Electrical and Electronic Engineering, Nanyang Technological University,
Singapore 639798 (email: csun002@e.ntu.edu.sg, gqhu@ntu.edu.sg).}}
\maketitle

\begin{abstract}
This paper studies the online convex optimization problem by using an Online Continuous-Time Nesterov Accelerated Gradient method  (OCT-NAG). We show that the continuous-time dynamics generated by the online version of the Bregman Lagrangian achieves a constant static regret $\frac{c}{\sigma}$ independent of $T$, provided that some boundedness assumptions on the objective functions and optimal solutions hold. To the best of the authors' knowledge, this is the lowest static regret in the literature (lower than $O(\text{log}(T))$). We further show that under the same assumptions, the dynamic regret of the algorithm is $O(T)$, which is comparable with the existing methods. Simulation results validate the effectiveness and efficiency of the method. Furthermore, the simulation shows that the algorithm performs well in terms of the dynamic regret for some specific scaling conditions. In addition, we consider the application of the proposed online optimization method in distributed online optimization problems, and show that the proposed algorithm achieves an $O(\sqrt{T})$ static regret, which is comparable with the existing distributed online optimization methods. Different from these methods, the proposed method requires neither the gradient boundedness assumption nor the compact constraint set assumption, which allows different objective functions and different optimization problems with those in the literature. A comparable dynamic regret is obtained. Simulation results show the effectiveness and efficiency of the distributed algorithm.

\end{abstract}


\thispagestyle{empty}

\section{Introduction}

Optimization plays an important role in many areas such as robotics, machine learning, finance, power electronics, and military, etc. Till now, numerous optimization methods have been developed with applications to different scenarios and for different algorithm requirements. Some popular methods include first-order gradient based methods (Gradient Descent, Projected Subgradient, Mirror Descent, Proximal Gradient, Alternating Direction Method of Multiplier (ADMM), etc), second-order Hessian based methods (Newton's method, Gauss-Newton, etc), Conjugate Gradient Method and stochastic optimization methods (Stochastic Gradient Descent, etc), to name a few \cite{boyd2004convex, bertsekas1997nonlinear}.

Convergence rate is one of the most important criteria to evaluate the performance of an optimization algorithm. Among all first-order gradient based methods, Nesterov Accelerated Gradient method has the optimal convergence rate according to first-order oracle complexity \cite{nesterov27method,nesterov2013introductory}. The acceleration feature makes the method become one of the most popular optimization methods and has a wide range of applications in many areas such as deep learning. Recently, there are also some further studies on Nesterov Accelerated Gradient methods in multi-agent systems. For example, \cite{qu2017accelerated} studied the applications of the algorithm in distributed optimization problems and \cite{tatarenko2018accelerated} applied the algorithm to distributively seek a Nash equilibrium of a noncooperative game.

The Nesterov Accelerated Gradient method is being studied and explained from a continuous-time dynamical system perspective \cite{su2014differential, pnas}. In \cite{su2014differential}, a second-order ordinary differential equation (ODE) was derived which has an approximate equivalence to the discrete-time Nesterov Accelerated Gradient method in terms of the convergence rate. Surprisingly, the continuous-time ODE generates a family of schemes with higher order convergence rates $O(\frac{1}{k^r}), r\geq 2$ than $O(\frac{1}{k^2})$ in \cite{nesterov27method}. In \cite{pnas}, a Bregman Lagrangian function was defined. Taking the Euler-Lagrange curve of the function generates a family of continuous-time ODEs, which are more general than the ODEs in \cite{su2014differential}. It was proven that these ODEs can achieve an arbitrary polynomial-order convergence rate and an exponential convergence rate under different scaling conditions (similar to the step-size concept) for a convex function. These interesting findings make the continuous-time methods attract more and more attentions recently. For example,
\cite{vassilis2018differential} explored the differential inclusion version of the algorithm by using nonsmooth analysis theory.

Online optimization was first defined in machine learning, which solves a class of optimization problems where the objective function is time-varying \cite{hazan2016introduction}. At each iteration $k$, the decision maker makes a decision $x_k$ according to the previous knowledge. After that, the decision maker receives its objective function $f_k$ of iteration $k$. In other words, the decision maker uses the previous information to help make a decision for the present optimization problem. For example, for an Online Gradient Descent method,
\begin{eqnarray} 
x(t_{k+1})=P_{\mathcal{C}}(x(t_{k})-\eta_k \nabla f_{t_k}(x(t_{k}))), 
\end{eqnarray}
where $x(t_{k})$ is the decision variable at iteration $k$, $P_{\mathcal{C}}$ is the projection onto the feasible set, $\eta_k $ is the step-size and $f_{t_k}(x(t_{k}))$ is the objective function at iteration $k$. The information at iteration $k$ is used to generate a decision for iteration $k+1$. The online optimization problem was also studied in a continuous-time sense. In \cite{popkov2005gradient}, a continuous-time gradient-based method was proposed as follows:
\begin{eqnarray} 
\dot{x}=-\gamma \nabla f_t(x), 
\end{eqnarray}
where $\gamma$ is a fixed gain. The present information $f_t(x)$ is used to help generate the guidance (i.e., $\dot{x}$) for the next moment. In addition, some distributed continuous-time online optimization algorithms were proposed in \cite{Chao, rahili2017distributed}. 

The regret is usually used to evaluate the effectiveness of an algorithm. The static regret refers to the accumulated difference between the cost function computed according to the value of the algorithm variable and the cost function computed according to the best fixed point that minimizes the accumulated cost function \cite{hazan2016introduction}, i.e., for a positive integer $K$ and $T=Kh$ with $h$ being the sampling period,
\begin{eqnarray}
\mathcal{R}_{ds}(T)=\sum_{t_k=T_0}^T (f_{t_k}(x(t_k))- f_{t_k}(\tilde{x}(T))), \label{regret0}
\end{eqnarray}
where $\tilde{x}(T)=\text{argmin}_{x\in\mathcal{C}}\sum_{k=0}^K f_{t_k}(x)$ is a comparator  in an offline sense. The dynamic regret is an online comparator which can be defined as 
\begin{eqnarray}
\mathcal{R}_{dd}(T)=\sum_{t_k=T_0}^T (f_{t_k}(x(t_k))- f_{t_k}(x^*(t_k))), \label{regret1}
\end{eqnarray}
where $x^*(t_k)$ is the optimal solution of the optimization problem at iteration $k$. Similarly, a recent paper defines the regrets in a continuous-time perspective \cite{lee2016distributed}. 

Some popular online optimization algorithms include Online Gradient Descent, Adaptive Gradient Method (AdaGrad), Online Mirror Descent, Online Newton Step, Follow-The-Approximate-Leader (FTAL), etc.  Under the assumption that the gradients are bounded, it was shown that the static regret for Online Gradient Descent is $O(\sqrt{T})$ \cite{zinkevich2003online} (under the gradient boundedness assumption, the bound is tight ), and the regret bound can be relaxed to $O(\text{log}(T))$ for strongly convex functions. Under the gradient boundedness assumption and the exp-concavity assumption on the convex cost functions (which is more general than strong convexity \cite{hazan2007logarithmic}), it was proven that the Online Newton Step and FTAL can achieve $O(\text{log}(T))$ static regrets. See Table I for the static regrets and conditions for these algorithms.

Distributed optimization aims to solve a class of network optimization problems with multiple agents where the objective function is the sum of all the agents' local objective functions \cite{nedic2010constrained, Chao, liu2017convergence}. The agents collaborate with each other to solve the optimization problem by using only neighboring information in the communication topology. Distributed online optimization deals with time-varying cost functions, and is attracting more and more attentions in many areas such as decentralized tracking \cite{shahrampour2017distributed} and swarm robots \cite{rahili2017distributed}. In addition to the above-mentioned continuous-time algorithms \cite{Chao, rahili2017distributed}, there are some discrete-time algorithms developed for distributed online optimization problems. For example, a distributed dual averaging method was proposed in \cite{hosseini2013online} and an $O(\sqrt{T})$ static regret was proven. A Distributed Online Mirror Descent method was proposed in \cite{shahrampour2017distributed} and the algorithm achieves an $O(\sqrt{T(1+V_T)})$ dynamic regret, where $V_T$ is a $T$-related function. An online subgradient descent algorithm with proportional-integral disagreement feedback was proposed in \cite{mateos2014distributed}. $O(\sqrt{T})$ and $O(\text{log}(T))$ regrets were proven under different assumptions. A saddle point algorithm was proposed in \cite{koppel2015saddle} with an $O(\sqrt{T})$ static regret. In 
\cite{hosseini2016online}, Distributed Weighted Dual Averaging was proposed with an $O(\sqrt{T})$ static regret.

Based on the above knowledge, this work aims to answer the following question: \textit{Whether the continuous-time Nesterov Accelerated Gradient method could be used to solve online optimization problems? Under what assumptions, what kind of static regrets and dynamic regrets could be guaranteed? What is the performance of the algorithm in distributed online optimization problems?}

Motivated by \cite{pnas}, We analyze the static regret and the dynamic regret of the Euler-Lagrange ODE generated by the online version of the Bregman Lagrangian function. Then, we design a distributed algorithm for distributed online optimization problems and analyze its performance. The main contributions are listed as follows:

(1) We prove that under some scaling conditions and under some boundedness assumptions, the static regret of the OCT-NAG is $\frac{c}{\sigma}$ which is a constant independent of $T$. This relaxes the static regrets in the literature. As a comparison, in addition to the $O(\text{log}(T))$ regrets for the algorithms in Table I, the most related work may be \cite{minimax} where a $(\text{log} (T)-\text{log} (\text{log} (T))+O(\frac{\text{log} (\text{log} (T))}{\text{log} (T)}))$ regret was found for some quadratic functions with some additional assumptions on the cost functions. To the best of the authors' knowledge, the proposed method achieves the lowest static regret bound in the literature.

(2) We prove that under the same scaling conditions, the dynamic regret of OCT-NAG is $O(T)$. This regret is comparable with the regrets of the existing methods in the literature. Most of the algorithms achieve similar dynamic regrets. For example, the Online Gradient Descent \cite{zinkevich2003online} achieves an $O(T)$ dynamic regret for fixed step-sizes. The Dynamic Mirror Descent in \cite{7044563} achieves an $O(\sqrt{T}(1+V_T))$ dynamic regret. A restarted Online Gradient Descent method in \cite{besbes} was proven to have a $O(T^{\frac{2}{3}}(V_T+1)^{\frac{1}{3}})$ dynamic regret.   Moreover, we show in the simulation that for some specific scaling conditions (i.e., $e^{\alpha_t+\beta_t}-\dot{\beta}_te^{\beta_t}=\sigma(t+b_0)$), the algorithm has a better performance in terms of the dynamic regrets than other compared algorithms.

(3) The proposed Distributed Online Continuous-Time Nesterov Accelerated Gradient method (DOCT-NAG) achieves a comparable static regret with the methods in the literature ($O(\sqrt{T})$). The algorithm requires neither the gradient boundedness assumption nor the projection onto a compact set, which could be used to solve different optimization problems with the literature \cite{hosseini2013online,shahrampour2017distributed,mateos2014distributed,koppel2015saddle,hosseini2016online}. 

(4) The discrete-time Nesterov Accelerated Gradient methods have been studied from various perspectives for static optimization/game problems \cite{qu2017accelerated, tatarenko2018accelerated}, and there are also some literature on their applications to online optimization\cite{NIPS2009_3817}. However, different from  \cite{qu2017accelerated, tatarenko2018accelerated, NIPS2009_3817}, we investigate the \textit{Online Continuous-Time} Nesterov Accelerated Gradient method and obtain lower static regrets than \cite{NIPS2009_3817}. In addition, most of the existing literature on online optimization are based on discrete-time algorithms. The continuous-time framework in this work allows for the possibility of solving optimization problems using well-established theories in ODE and control theory.

The rest of the paper is organized as follows: In Section \ref{S2}, notations and preliminary knowledge on graph theory and convex optimization is given. In Section \ref{S3}, the static regret and the dynamic regret of OCT-NAG are analyzed. In Section \ref{dis}, the static regret and the dynamic regret of DOCT-NAG are given. In Section \ref{S4}, simulation results show the effectiveness and efficiency of the proposed algorithms. Finally, Section \ref{S5} concludes the paper.

\section{Notations and Preliminaries} \label{S2}    

\textbf{Notations}: Letting $\mathcal{C}$ being a convex set, the notation $\mathop{\text{argmin}}\limits_{x\in\mathcal{C}} f(x)$ represents the decision variable value $x\in\mathcal{C}$ such that the function $f(x)$ takes its minimum. $\left\Vert \cdot\right\Vert$ is the Euclidean norm and $\left\vert\cdot\right\vert$ is the absolute value. log$(\cdot)$ represents the natural logarithm. ``$0$" represents a zero scalar or vector with a appropriate dimension. $1_N$ is an $N\times 1$ column vector with all elements being 1.

\textbf{Graph Theory}:  Let $\mathcal{G}=\left\{  \mathcal{V},\mathcal{E}\right\}  $ denote an
undirected graph, where $\mathcal{V=}\left\{  1,...,N\right\}  $ is the
vertex set and $\mathcal{E\subset V\times V}$ is the edge set.
$\mathcal{N}_{i}=\left\{  j\in\mathcal{V\mid}(j,i)\in\mathcal{E}\right\}  $
denotes the set of neighbors of vertex $i$. A path is referred by the sequence
of its vertices. Path $\mathcal{P}$ between $v_{0}$ and $v_{m}$\ is the
sequence $\left\{  v_{0},...,v_{m}\right\}  $ where $(v_{i-1},v_{i}%
)\in\mathcal{E}$ for $i=1,...,m$ and the vertices are distinct. The number $m$
is defined as the length of path $\mathcal{P}$. Graph $\mathcal{G}$\ is
connected if for any two vertices, there is a path in $\mathcal{G}$. A matrix
$A=\left[  a_{ij}\right]  \in%
%TCIMACRO{\U{211d} }%
%BeginExpansion
\mathbb{R}
%EndExpansion
^{N\times N}$ denotes the adjacency matrix of $\mathcal{G}$, where $a_{ij}>0$
if and only if $(j,i)\in\mathcal{E}$ else $a_{ij}=0$. In this paper, we
suppose that $a_{ii}=0.$  The degree of a vertex $i$ in an undirected graph is the sum of the weights in $\mathcal{N}_{i}$. A matrix $L\triangleq D-A\in%
%TCIMACRO{\U{211d} }%
%BeginExpansion
\mathbb{R}
%EndExpansion
^{N\times N}$ is called the Laplacian matrix of $\mathcal{G}$, where
$D=\text{diag} { \{d_{ii}\}} \in%
%TCIMACRO{\U{211d} }%
%BeginExpansion
\mathbb{R}
%EndExpansion
^{N\times N}$ is a diagonal matrix with $d_{ii}=\sum\nolimits_{j=1}^{N}a_{ij}
$ \cite{RenTAC05}.

\textbf{Convex Optimization}:  For a continuously differentiable function $f(x)$, it is convex over a convex set $\mathcal{C}$ if and only if $f(x)-f(y)\geq (\nabla_x f(y))^T (x-y)$ for all $x,y\in\mathcal{C}$.

\section{Online Continuous-time Nesterov Accelerated Gradient} \label{S3}  

\subsection{Problem Formulation and Regrets} \label{S3.1}

Consider a continuous-time online optimization problem:
\begin{eqnarray}
\text{min}_{x\in\mathcal{C}} f_t(x), \label{problem1}
\end{eqnarray}
where $\mathcal{C}\subset \mathbb{R}$ is a convex set,  $x\in\mathbb{R}$ is the decision variable, $f_t(x)$ is a time-varying, continuously differentiable cost function, which is convex at any time $t$, and smooth in $t$. The objective is to design an algorithm to obtain and track the optimal solution of the optimization problem at each time instant. We assume that at any time $t$, the function $f_t(x)$ has at least one minimizer.

The static regret function of the optimization problem at time $t=T>T_0$ can be defined as
\cite{lee2016distributed}:
\begin{eqnarray}
\mathcal{R}_s(T)=\int_{T_0}^T (f_t(x(t))- f_t(\tilde{x}(T)))dt, \label{regret}
\end{eqnarray}
where $\tilde{x}(T)$ is a constant with respect to $T$ defined as 
\begin{eqnarray}
\tilde{x}(T)=\mathop{\text{argmin}}\limits_{x\in\mathcal{C}}\int_{T_0}^T f_t(x)dt, \label{offline}
\end{eqnarray}
or equivalently, 
\begin{eqnarray}
\int_{T_0}^T f_t(\tilde{x}(T))dt=\text{min}_{x\in\mathcal{C}} \int_{T_0}^T f_t(x)dt.
\end{eqnarray}

In (\ref{regret}), the first term $\int_{T_0}^Tf_t(x(t))dt $ represents the integration of the cost function computed according to the value of the algorithm variable from $T_0$ to $T$, and the second term $\int_{T_0}^T f_t(\tilde{x}(T))dt$ represents the integration of the cost function computed according to the best fixed point that minimizes the integrated cost function from $T_0$ to $T$ (an offline minimizer).

The dynamic regret function of the optimization problem at time $t=T$ can be defined as
\begin{eqnarray}
\mathcal{R}_d(T)=\int_{T_0}^T (f_t(x(t))- f_t(x_t^*))dt, \label{regret2}
\end{eqnarray}
where $x_t^*$ is a time-varying optimal solution defined as 
\begin{eqnarray}
x_t^*=\mathop{\text{argmin}}\limits_{x\in\mathcal{C}} f_t(x), \label{eq10}
\end{eqnarray}
or equivalently, 
\begin{eqnarray}
f_t(x_t^*)=\text{min}_{x\in\mathcal{C}} f_t(x).
\end{eqnarray}

In this section, we focus on the case $\mathcal{C}=\mathbb{R}$, i.e., the optimization problem defined
in  (\ref{problem1}) is an unconstrained online optimization problem.

\subsection{Algorithm Design}

Motivated by \cite{pnas}, consider the following Euler-Lagrange equation 
\begin{eqnarray}
\ddot{x}+(e^{\alpha_t}-\dot{\alpha}_t)\dot{x}+e^{2\alpha_t+\beta_t}\nabla f_t(x)=0, \label{algorithm}
\end{eqnarray}
where $\nabla f_t(x)$ is the time-varying gradient of the time-varying cost function $f_t(x)$ at time $t$, $\alpha_t$ and $\beta_t$ are time-varying parameters determined later. For simplifying the notation and without losing generality, we consider the Bregman Lagrangian \cite{pnas} in an Euclidean setting as in (\ref{algorithm}). In addition, in (\ref{algorithm}), we fix the Lagrangian damping to satisfy the ideal scaling condition. See \cite{pnas} on how to obtain the differential equation described in (\ref{algorithm}). The contribution of this section compared with \cite{pnas} is the extension of the equation to an online version, and the following regret analysis.

\subsection{Static Regret} \label{S3_1}

To analyze the static regret of (\ref{algorithm}), let $z\in\mathbb{R}$ be any constant. Define a function $V_z(t,x,\dot{x})\in\mathbb{R}$ as
\begin{eqnarray}
V_z=\frac{1}{2}(x+e^{{-\alpha}_t} \dot{x}-z)^2
+e^{\beta_t}(f_t(x)-f_t(z)).
\end{eqnarray}

Note that $V_z(t,x,\dot{x})$ is not necessarily larger than or equal to $0$. Furthermore,
\begin{eqnarray}
V_{z}&=&\frac{1}{2}(x+e^{{-\alpha}_t} \dot{x}-z)^2 \notag \\ & &+e^{\beta_t}(f_t(x)-f_t(x^*_t)) \notag \\
& &+e^{\beta_t}(f_t(x^*_t)-f_t(z)) \notag \\
&\geq& e^{\beta_t}(f_t(x^*_t)-f_t(z)), \label{Lya}
\end{eqnarray}
where $x^*_t$ was defined in (\ref{eq10}).

Taking the derivative of $V_z$ gives
\begin{eqnarray}
\dot{V}_z&=&(x+e^{{-\alpha}_t} \dot{x}-z)(\dot{x}-\dot{\alpha}_te^{{-\alpha}_t}\dot{x}+e^{{-\alpha}_t} \ddot{x}) \notag \\
& &+\dot{\beta}_te^{\beta_t}(f_t(x)-f_t(z))+e^{\beta_t}(\dot{x}\nabla f_t(x)
 \notag \\
& &+\nabla_t f_t(x)-\nabla_t f_t(z)) \notag \\
&=&(x+e^{{-\alpha}_t} \dot{x}-z)(-e^{\alpha_t+\beta_t}\nabla f_t(x)) \notag \\
& &+\dot{\beta}_te^{\beta_t}(f_t(x)-f_t(z))+e^{\beta_t}\dot{x}\nabla f_t(x)  \notag \\
& &+e^{\beta_t} (\nabla_t f_t(x)-\nabla_t f_t(z))  \notag \\
&\leq&-(e^{\alpha_t+\beta_t}-\dot{\beta}_te^{\beta_t})(f_t(x)-f_t(z))  \notag \\
& &+e^{\beta_t}(\nabla_t f_t(x)-\nabla_t f_t(z)), \label{derivative}
\end{eqnarray} 
where in the last step we use the fact that for any $x$ and $y$, $f_t(x)-f_t(y)\leq \nabla f_t(x)(x-y)$ for a convex function $f_t$.

Integrating both sides of (\ref{derivative}) from $T_0$ to $T>T_0$ gives
\begin{eqnarray}
& & \int_{T_0}^T(e^{\alpha_t+\beta_t}-\dot{\beta}_te^{\beta_t})(f_t(x(t))-f_t(z)) dt  \notag \\
&\leq& V_z(T_0,x(T_0),\dot{x}(T_0))-V_z(T,x(T),\dot{x}(T)) \notag \\
& &+\int_{T_0}^T e^{\beta_t}(\nabla_t f_t(x(t))-\nabla_t f_t(z))dt. \label{equation}
\end{eqnarray} 

\begin{remark}
	For (\ref{equation}), let $f_t(x)\equiv f(x)$ which is a static objective function, and let $z=x^*$ where $f(x^*)$ is the minimal function value of $f(x)$. Then, $V_z(T)\leq V_z(T_0)$ if $e^{\alpha_t}-\dot{\beta}_t\geq 0$. This conclusion is consistent with Theorem 1.1 of \cite{pnas}.
\end{remark}

According to (\ref{Lya}), we have 
\begin{eqnarray}
V_{z}(T,x(T),\dot{x}(T))\geq e^{\beta_T}(f_T(x^*_T)-f_T(z)). \label{equation3}
\end{eqnarray} 

In this section, we study the static regret defined in (\ref{regret}). To this end, we consider the following scaling condition: 
\begin{eqnarray}
\sigma=(e^{\alpha_t+\beta_t}-\dot{\beta}_te^{\beta_t}), \label{f}
\end{eqnarray} 
with $\sigma$ being a positive constant \footnote{Instead, one can also define $e^{\alpha_t+\beta_t}-\dot{\beta}_te^{\beta_t}=g(t)=\sigma t$ to obtain the regret in an average sense. Note that in this work, we only consider some specific scaling conditions. There may exist other scaling conditions that satisfy the conclusions in this work.}. 

To facilitate the following analysis, we make the following assumptions:
\begin{assumption}
	$x^*_t$ exists and is bounded for $t\in[T_0,T]$ with the bound independent of $T$, i.e., there exists a positive constant $c_0$ independent of $T$ such that $\left\vert x^*_t\right \vert \leq c_0$. \label{a1}
\end{assumption}

\begin{assumption}
 $\tilde{x}(T)= \mathop{\text{argmin}}\limits_{x\in\mathbb{R}} \int_{T_0}^T f_t(x) dt $ exists and is bounded with the bound independent of $T$.  \label{a2}
\end{assumption}

\begin{assumption}
	If $ x$ is bounded, $ f_t(x) $ is bounded for $t\in[T_0,T]$ with the bound independent of $T$. \label{a3}
\end{assumption}

\begin{assumption}
	$\nabla_t f_t(x)$ is bounded for $t\in[T_0,T]$ with the bound independent of $T$ and $x$. \label{boundedness}
\end{assumption}
  
\begin{remark} 
 (1) Assumption \ref{a1} is a standard assumption which implies the boundedness of the optimal solutions and the feasibility of the optimization problem. In a static optimization problem, it is usually assumed that a finite optimal solution exists. Assumption \ref{a1} is an extended form of this assumption in an online optimization problem.
  
  (2) Assumption \ref{a2} implies that the optimal solutions of the integration function of $f_t$ are also bounded, i.e., the offline minimizers exist at any time and are uniformly bounded. This assumption is mild. In fact, $\tilde{x}(T)$ serves only as a comparator of the algorithm. If this assumption does not hold, the offline minimizer could be infinity when $t$ tends to infinity. In this case, the offline optimization problem may not have a finite optimal solution. This implies that the static regret is meaningless for algorithm comparisons.
  
  (3) Assumption \ref{a3} implies that the boundedness of the time-related component in the objective function. For example, the component could be $\text{sin}(t)$ or $\text{tanh}(t)$.
  
  (4) Assumption \ref{boundedness} restricts the partial derivative of the function with respect to $t$ to be uniformly bounded for $x$ and $t$. Note that Assumption \ref{boundedness} is similar to Assumption 9.1 of \cite{simonetto2017time} on discrete-time online primal methods, both of which imply that the variation of the cost function with respect to $t$ cannot be too fast.
  
    (5) A typical example that satisfies Assumptions \ref{a1}-\ref{boundedness} is $f_t(x)=x^2+\text{sin}(t)\text{sin}(x)$. See Section \ref{S4} for more details on this point.
\end{remark}

\begin{remark} 
Assumptions \ref{a3} and \ref{boundedness} are restrictive compared with other algorithms, e.g., Online Gradient Descent. However, in the following, we show that for objective functions satisfying Assumptions \ref{a1}-\ref{boundedness}, the algorithm has a lower regret. Thus, the algorithm is of its own value despite of its restrictive assumptions. In addition, Assumptions \ref{a3} and \ref{boundedness} are related to the boundedness of some functions/components in the decision variable and $t$. It is interesting to study the algorithm when the feasible set $\mathcal{C}$ is a compact set.  
\end{remark}

Note that for any fixed $T$, $\tilde{x}(T)$ is fixed although it is related to $T$ (its bound is independent of $T$). According to (\ref{equation}), (\ref{equation3}), and (\ref{f}), and letting $z=\tilde{x}(T)$ gives
\begin{eqnarray}
& & \sigma\int_{T_0}^T(f_t(x(t))-f_t(\tilde{x}(T))) dt  \notag \\
&\leq& V_{\tilde{x}}(T_0,x(T_0),\dot{x}(T_0))-e^{\beta_T}(f_T(x^*_T)-f_T(\tilde{x}(T))) \notag \\
& &+\int_{T_0}^T e^{\beta_t}(\nabla_t f_t(x(t))-\nabla_t f_t(\tilde{x}(T)))dt. \label{equation2}
\end{eqnarray} 

To obtain a polynomial convergence rate, let 
\begin{eqnarray}
e^{\beta_t}=(t+b_0)^{m}, \label{eq}
\end{eqnarray}
 where $m<0$ is a constant which influences the convergence rate and $b_0>(\frac{-m}{\sigma})^{\frac{1}{1-m}}$ is a positive constant that is used to adjust the sign of some variables. Then, it can be calculated that
\begin{eqnarray}
\dot{\beta}_t&=&m(t+b_0)^{-1}, \notag \\
e^{\alpha_t}&=&m(t+b_0)^{-1}+\sigma(t+b_0)^{-m}>0, \notag \\
\dot{\alpha_t}&=&(m(t+b_0)^{-1}+\sigma(t+b_0)^{-m})^{-1} \notag \\
& &\times (-m(t+b_0)^{-2}-m\sigma(t+b_0)^{-m-1}). \label{fff}
\end{eqnarray}

Based on Assumptions \ref{a1}-\ref{a3}, there exists a positive constant $c_1$ such that 
\begin{eqnarray}
\left\vert f_t(x^*_t)-f_t(\tilde{x}(t))\right\vert \leq c_1  \label{a4}
\end{eqnarray}
hold for any $t\in[T_0,\infty)$.

According to (\ref{equation2}), it can be obtained that
\begin{eqnarray}
& & \int_{T_0}^T(f_t(x(t))-f_t(\tilde{x}(T))) dt  \notag \\
&\leq& \frac{1}{\sigma}\left\vert V_{\tilde{x}}(T_0,x(T_0),\dot{x}(T_0))\right\vert+ \frac{c_1}{\sigma}(T+b_0)^{m} \notag \\
& &+ \frac{1}{\sigma}\int_{T_0}^T e^{\beta_t}(\nabla_t f_t(x(t))-\nabla_t f_t(\tilde{x}(T)))dt. \label{equation4}
\end{eqnarray}

Thus, if $-1<m<0$ or $m< -1$,
\begin{eqnarray}
& & \int_{T_0}^T(f_t(x(t))-f_t(\tilde{x}(T))) dt  \notag \\
&\leq& \frac{1}{\sigma}\left\vert V_{\tilde{x}}(T_0,x(T_0),\dot{x}(T_0))\right\vert+\frac{c_1}{\sigma}(T+b_0)^{m} \notag \\
& &+\frac{c_2}{\sigma|m+1|}(|(T+b_0)^{m+1}|+|(T_0+b_0)^{m+1}|), \label{cequation1}
\end{eqnarray} 
where $c_2$ is some positive constant satisfying $\left\vert\nabla_t f_t(x(t))-\nabla_t f_t(\tilde{x}(T))\right\vert\leq c_2$, and the boundedness can be guaranteed by Assumption \ref{boundedness}.

If $m=-1$,
\begin{eqnarray}
& & \int_{T_0}^T(f_t(x(t))-f_t(\tilde{x}(T))) dt  \notag \\
&\leq& \frac{1}{\sigma}\left\vert V_{\tilde{x}}(T_0,x(T_0),\dot{x}(T_0))\right\vert+\frac{c_1}{\sigma}(T+b_0)^{m} \notag \\
& &+\frac{c_2}{\sigma}\text{log}\frac{T+b_0}{T_0+b_0}. \label{cequation2}
\end{eqnarray} 

Based on Assumptions \ref{a1}-\ref{boundedness}, $V_{\tilde{x}}(T_0,$ $x(T_0),$ $\dot{x}(T_0))$ can be bounded by a constant $c_3>0$ independent of $T$.

Then, we have the following conclusion:

\begin{theorem}
	Assume that Assumptions \ref{a1}-\ref{boundedness} hold. The Online Continuous-Time Nesterov Accelerated Gradient dynamics in
	(\ref{algorithm}) satisfies that for all $T>T_0$,
	\begin{eqnarray}
	\int_{T_0}^T(f_t(x(t))-f_t(\tilde{x}(T))) dt\leq \frac{c}{\sigma}, \label{theoremcondition}
	\end{eqnarray}
	with $c$ being some positive constant independent of $T$ and $\sigma$ being any positive parameter,
	provided that the conditions in  (\ref{f}), (\ref{eq}), and (\ref{fff}) hold and $m<-1$. \label{theorem1}
\end{theorem}

\begin{remark}
	The conclusion in Theorem \ref{theorem1} does not imply that $\lim_{t\rightarrow \infty}(f_t(x(t))-f_t(\tilde{x}(T)))=0$. For example,
	if for all $T_0\leq t \leq T$, $f_t(x(t))-f_t(\tilde{x}(T))\leq 0$, then, (\ref{theoremcondition}) holds. In this case, it implies that the algorithm will always obtain a smaller function value than the offline minimizer.
	\end{remark}

\subsection{Dynamic Regret} \label{aaa}

Considering the scaling condition  (\ref{f}), based on (\ref{equation2}), 
\begin{eqnarray}
& & \sigma\int_{T_0}^T(f_t(x(t))-f_t({x}_t^*)) dt  \notag \\
&\leq& -\sigma\int_{T_0}^T(f_t({x}_t^*)-f_t(\tilde{x}(T))) dt +V_{\tilde{x}}(T_0,x(T_0),\dot{x}(T_0)) \notag \\
& &-e^{\beta_T}(f_T(x^*_T)-f_T(\tilde{x}(T)))+\int_{T_0}^T e^{\beta_t}(\nabla_t f_t(x(t)) \notag \\
& &-\nabla_t f_t(\tilde{x}(T)))dt. \label{equation2d}
\end{eqnarray} 

Similar to Section \ref{S3_1}, letting (\ref{eq}) hold, it can be obtained that if $-1<m<0$ and $m<-1$,
\begin{eqnarray}
& & \int_{T_0}^T(f_t(x(t))-f_t({x}_t^*)) dt  \notag \\
&\leq&(T-T_0)c_1+\frac{1}{\sigma}\left\vert V_{\tilde{x}}(T_0,x(T_0),\dot{x}(T_0))\right\vert+\frac{c_1}{\sigma}(T+b_0)^{m} \notag \\
& &+\frac{c_2}{\sigma|m+1|}((T+b_0)^{m+1}+(T_0+b_0)^{m+1}), \notag 
\end{eqnarray} 
and if  $m=-1$,
\begin{eqnarray}
& & \int_{T_0}^T(f_t(x(t))-f_t({x}_t^*)) dt  \notag \\
&\leq& (T-T_0)c_1+\frac{1}{\sigma}\left\vert V_{\tilde{x}}(T_0,x(T_0),\dot{x}(T_0))\right\vert+\frac{c_1}{\sigma}(T+b_0)^{m} \notag \\
& &+\frac{c_2}{\sigma}\text{log}\frac{T+b_0}{T_0+b_0}. \notag 
\end{eqnarray} 

\begin{theorem}
	Assume that Assumptions \ref{a1}-\ref{boundedness} hold. The OCT-NAG algorithm  in
	(\ref{algorithm}) satisfies that for all $T>T_0$,
	\begin{eqnarray}
	\int_{T_0}^T(f_t(x(t))-f_t(x_t^*)) dt\leq (T-T_0)c_1+\frac{c}{\sigma}, \label{theoremcondition2}
	\end{eqnarray}
	provided that the conditions in  (\ref{f}), (\ref{eq}), and (\ref{fff}) hold and $m< -1$. \label{theorem2}
\end{theorem}

\begin{table*}
	\centering
	\caption{Static Regret Bounds and Conditions of Some Online Optimization Algorithms}
\scalebox{0.8}{
\begin{tabular}{c|c|c|c}
	\hline
	Algorithms & References & Bounds & Conditions \\
	\hline  
	Online Gradient Descent & \cite{hazan2016introduction} & $\sum_{t_k=T_0}^T (f_{t_k}(x(t_k))- f_{t_k}(\tilde{x}(T)))\leq O(\sqrt{T})$ & Convex, Gradient Boundedness\\
	\hline 
	Online Gradient Descent & \cite{hazan2016introduction} & $\sum_{t_k=T_0}^T (f_{t_k}(x(t_k))- f_{t_k}(\tilde{x}(T)))\leq O(\text{log}{T})$ & Strongly Convex, Gradient Boundedness\\
	\hline 
	AdaGrad & \cite{hazan2016introduction, duchi2011adaptive} & $\sum_{t_k=T_0}^T (f_{t_k}(x(t_k))- f_{t_k}(\tilde{x}(T)))\leq O(\sqrt{T})$ & Convex, Gradient Boundedness\\
	\hline 
	Online Mirror Descent & \cite{hazan2016introduction, shalev2012online} & $\sum_{t_k=T_0}^T (f_{t_k}(x(t_k))- f_{t_k}(\tilde{x}(T)))\leq O(\sqrt{{T}})$ & Convex, Gradient Boundedness\\
	\hline 
	Online Newton Step &  \cite{hazan2006efficient, hazan2007logarithmic} & $\sum_{t_k=T_0}^T (f_{t_k}(x(t_k))- f_{t_k}(\tilde{x}(T)))\leq O(\text{log}T)$ & $e^{-\alpha f_{t_k}(x)}$ is concave, Gradient Boundedness\\
	\hline 
	FTAL &  \cite{hazan2007logarithmic} & $\sum_{t_k=T_0}^T (f_{t_k}(x(t_k))- f_{t_k}(\tilde{x}(T)))\leq O(\text{log}T)$ & $e^{-\alpha f_{t_k}(x)}$ is concave, Gradient Boundedness\\
	\hline 
	Exponentially-Weighted-Online-Optimization &  \cite{hazan2007logarithmic} & $\sum_{t_k=T_0}^T (f_{t_k}(x(t_k))- f_{t_k}(\tilde{x}(T)))\leq O(\text{log}T)$ & $e^{-\alpha f_{t_k}(x)}$ is concave\\
	\hline 
	This work &  & $\int_{t=T_0}^T (f_{t}(x(t))- f_t(\tilde{x}(T)))dt\leq \text{positive constant} $ & Convex, Boundedness Assumptions\\
	\hline
\end{tabular}
}
\label{table1}
\end{table*}

Consider the following scaling condition
\begin{eqnarray}
e^{\alpha_t+\beta_t}-\dot{\beta}_te^{\beta_t}=\sigma(t+b_0)^p, \label{f_2}
\end{eqnarray} 
with $p\geq 1$ being a positive constant, and $b_0\geq (\frac{-m}{\sigma})^{\frac{1}{1-m+p}}$.

Letting (\ref{eq}) hold, it can be obtained that
\begin{eqnarray}
\dot{\beta}_t&=&m(t+b_0)^{-1}, \notag \\
e^{\alpha_t}&=&m(t+b_0)^{-1}+\sigma(t+b_0)^{-m+p}>0, \notag \\
\dot{\alpha_t}&=&(m(t+b_0)^{-1}+\sigma(t+b_0)^{-m+p})^{-1} \notag \\
& &\times (-m(t+b_0)^{-2}-(m-p)\sigma(t+b_0)^{-m+p-1}). \notag \\
 \label{fff_2}
\end{eqnarray}

Similar to the above analysis, we have the following conclusion:

\begin{theorem}
	Assume that Assumptions \ref{a1}-\ref{boundedness} hold. The OCT-NAG algorithm in
	(\ref{algorithm}) satisfies that for all $T>T_0$,
	\begin{eqnarray}
	\int_{T_0}^T(t+b_0)^p(f_t(x(t))-f_t(x_t^*)) dt \notag \\ \leq((T+b_0)^{p+1}-(T_0+b_0)^{p+1})\frac{c_1}{p+1}+\frac{c_4}{\sigma} \label{theoremcondition1}
	\end{eqnarray}
	with $c_2$ being some constant, provided that the conditions in  (\ref{eq}), (\ref{f_2}),  and (\ref{fff_2}) hold and $m<-1$. \label{theorem3}
\end{theorem}

\begin{remark}
	The parameters selected in Theorem \ref{theorem3} does not provide a better theoretical regret bound than the parameters in Theorem \ref{theorem2}. However, it is found in the simulation that the algorithm with parameters in Theorem \ref{theorem3} may have a better performance. 
	\end{remark}

\section{Distributed Online Continuous-time Nesterov Accelerated Gradient Method} \label{dis}
\subsection{Problem Formulation and Regrets}

Consider that a group of agents collaborate with each other to solve the following distributed continuous-time online optimization problem:
\begin{eqnarray}
\text{min}_{x\in\mathbb{R}^n} f_t(x)=\sum_{i=1}^N f_{i,t}(x), \label{problem1d}
\end{eqnarray}
where $x\in\mathbb{R}^n$ is the decision variable, $f_{i,t}(x)$ is a time-varying, continuously differentiable local cost function, which is convex at any time $t$, and smooth in $t$, and $f_{t}(x)$ is the global cost function. 

Suppose that each agent can only get the information from its neighbors via an undirected and connected graph $\mathcal{G}$.  The objective is to design an algorithm to obtain and track the optimal solution of the optimization problem at each time instant in a distributed way. We assume that at any time $t$, the functions $f_t(x)$ and $f_{i,t}(x)$ have at least one minimizer.

Motivated by \cite{Gharesifard14}, let $\mathbf{x}_i\in\mathbb{R}^n$ be agent $i$'s estimation on the optimal solution. Then, the problem in (\ref{problem1d}) can be reformulated as follows:
\begin{eqnarray}
\text{min}_{\mathbf{x}_i\in\mathbb{R}^n} &F_t(\mathbf{x})=\sum_{i=1}^N f_{i,t}(\mathbf{x}_i), \notag \\
\text{subject to } &\tilde L \mathbf{x}=0,   \label{problem3}
\end{eqnarray}
where $\mathbf{x}=[\mathbf{x}_1^T,\cdots,\mathbf{x}_N^T]^T\in\mathbb{R}^{Nn}$,  $\tilde L\triangleq L\otimes I_n$ and $L$ is the Laplacian matrix of graph $\mathcal{G}$. 

Let $\mathbf{x}^*(t)$ be an optimal solution of (\ref{problem3}), it can be verified that $\mathbf{x}^*(t)=1_N \otimes x^*_t$ with $x^*_t=[x^*_{1,t},\cdots,x^*_{n,t}]^T\in\mathbb{R}^n$ being an optimal solution of (\ref{problem1d}).

To evaluate the performance of the distributed online optimization algorithm, we adopt the regrets defined in \cite{shahrampour2017distributed}, where the static regret was defined by
\begin{eqnarray}
\mathcal{R}_{s,dis}(T)=\int_{T_0}^T (\frac{1}{N}\sum_{j=1}^N\sum_{i=1}^N f_{i,t}(\mathbf{x}_j(t))- F_t(\mathbf{\tilde{x}}(T)))dt, \label{regretd1}
\end{eqnarray}
where $\mathbf{\tilde{x}}(T)=1_{N}\otimes\tilde{x}(T)$, $\tilde{x}(T)\in\mathbb{R}^n$ is defined by
\begin{eqnarray}
\tilde{x}(T)=\mathop{\text{argmin}}\limits_{x\in\mathbb{R}^n}\int_{T_0}^T f_t(x)dt, \label{offline1}
\end{eqnarray}
and the sum operator over the time is replaced by the integration operator for the continuous-time algorithm. 

The dynamic regret was defined by
\begin{eqnarray}
\mathcal{R}_{d,dis}(T)=\int_{T_0}^T (\frac{1}{N}\sum_{j=1}^N\sum_{i=1}^N f_{i,t}(\mathbf{x}_j(t))- F_t(\mathbf{x}^*(t)))dt. \label{regretd2}
\end{eqnarray}

\subsection{Algorithm Design}

The distributed updating law for agent $i$ is designed as 
\begin{align}
&  \ddot{\mathbf{x}}_i+(2e^{\alpha_t}-\dot{\alpha}_t)\dot{\mathbf{x}}_i+e^{2\alpha_t+\beta_t}\nabla f_{i,t}(\mathbf{x}_i) \notag \\
& +\sum_{j=1}^N a_{ij} k_1 e^{2\alpha_t}  (\mathbf{x}_i-\mathbf{x}_j)=0, \label{algorithm2}
\end{align}
where $k_1$ is a positive constant.

Compared with the centralized algorithm in (\ref{algorithm}), the term $\sum_{j=1}^N a_{ij} k_1 e^{2\alpha_t}  (\mathbf{x}_i-\mathbf{x}_j)$ is added for the equality constraint in (\ref{problem3}). In addition, the time-varying parameter for $\dot{\mathbf{x}}_i$ is changed to $2e^{\alpha_t}-\dot{\alpha}_t$ due to the requirement of the $T$-related boundedness of $\mathbf{x}_i$ in the following derivations.

Then, the concatenated form of (\ref{algorithm2}) can be written as 
\begin{align}
&\ddot{\mathbf{x}}+(2e^{\alpha_t}-\dot{\alpha}_t)\dot{\mathbf{x}}+e^{2\alpha_t+\beta_t}\nabla F_{t}(\mathbf{x}) \notag \\
&  + k_1 e^{2\alpha_t} \tilde L \mathbf{x}=0.   \label{algorithm3}
\end{align}

In this section, the following assumptions will be used.

\begin{assumption}
	$x^*_t$ exists and is bounded for $t\in[T_0,T]$ with the bound independent of $T$, i.e., there exists a positive constant $c_0'$ independent of $T$ such that $\left\Vert x^*_t\right \Vert \leq c_0'$. \label{a1dis}
\end{assumption}

\begin{assumption}
	$\tilde{x}(T)= \mathop{\text{argmin}}\limits_{x\in\mathbb{R}^n} \int_{T_0}^T f_t(x) dt $ exists and is bounded with the bound independent of $T$.  \label{a2dis}
\end{assumption}

\begin{assumption}
	$\mathbf{x}_{i}^\dagger(t)$ exists which is an optimal solution of $f_{i,t}(x)$ and is bounded for $t\in[T_0,T]$ with the bound independent of $T$, i.e., there exists a positive constant $c^\dagger$ independent of $T$ such that $\left\Vert \mathbf{x}_{i}^\dagger(t)\right \Vert \leq c^\dagger$. \label{as}
\end{assumption}

\begin{assumption}
	If $ x$ is bounded, $ f_{i,t}(x) $ is bounded for $t\in[T_0,T]$ with the bound independent of $T$, for all $i\in\mathcal{V}$. \label{a32}
\end{assumption}

\begin{assumption}
	$\nabla_t F_t(\mathbf{x})$ is bounded for $t\in[T_0,T]$ with the bound independent of $T$ and $\mathbf{x}$. \label{boundedness1}
\end{assumption}

%\begin{assumption} \cite{paternain2019distributed}
%	There exists a constant $\varrho>0$ independent of $T$ such that for all $t\in (T_0, T]$, 
%	\begin{eqnarray}
%	f_t(\tilde{x}(T))-\text{min}_{\mathbf{x}\in\mathbb{R}^{Nn}}F_t(\mathbf{x})\leq \varrho.
%	\end{eqnarray} \label{boundedness2}
%\end{assumption}
\begin{assumption} 
If for any $T_0<t\leq T$,  $\left\Vert x(t) \right\Vert \leq a_1\sqrt{t}+b_1$, then there exists a constant $m_0$ such that $\left\Vert \nabla f_{i,t}(x(t))\right\Vert \leq a_2 t^{m_0}+b_2$ for all $i\in\mathcal{V}$, where $a_1,a_2, b_1$, and $b_2$ are positive constants independent of $T$. \label{boundednessgra}
\end{assumption}

\begin{remark}
	Assumptions \ref{a1dis} and \ref{as} are similar to Assumption \ref{a1}.  Assumption  \ref{a2dis} is similar to Assumption  \ref{a2}. Assumption  \ref{a32} is similar to Assumption \ref{a3}. Assumption \ref{boundedness1} is similar to Assumption  \ref{boundedness}.
 Assumption \ref{boundednessgra} is a relaxed version of the gradient boundedness assumption \cite{mateos2014distributed}. In the existing literature, most of the work requires the gradient boundedness (or else, a compact feasible solution set is required and projection onto the compact set at each iteration is needed). Assumptions \ref{a1dis}-\ref{boundednessgra} are indeed not weaker than those in the existing work. However, the assumptions in this work allow some different objective functions with those in the literature.  An example that satisfies Assumptions \ref{a1dis}-\ref{boundednessgra} but doesn't satisfy the gradient boundedness assumption is $f_t(x)=\sum_{i=1}^N (i\times x_i^2+\text{sin}(0.1\times i\times t)\text{sin}(x_i))$.
\end{remark}

\subsection{Static Regret} \label{S4_1}

%Let $\mathbf{x}_e, \mathbf{y}_e$ be the quasi-steady state of (\ref{algorithm3}) satisfying $\ddot{\mathbf{x}}_e=0$, $\dot{\mathbf{x}}_e=0$, and $\dot{\mathbf{y}}_e=0$. Then, it gives that 
%\begin{align}
%e^{2\alpha_t+\beta_t}\nabla F_{t}(\mathbf{x}_e)+k_2 \mathbf{y}_e=0, \tilde L \mathbf{x}_e=0. \label{equi}
%\end{align}

%According to (\ref{equi}), there exists a constant vector $x_e^*\in\mathbb{R}^n$ such that $\mathbf{x}_e=1_N \otimes x_e^*$. Since $\sum_{i=1}^N \dot{{\mathbf{y}}}_i=0$, $\sum_{i=1}^N {\mathbf{y}}_i=\sum_{i=1}^N {\mathbf{y}}_i(0)=0$. Based on (\ref{equi}), it gives that $\sum_{i=1}^N \nabla f_{i,t}(\mathbf{x}_e) =0$. Then, based on the convexity, $\mathbf{x}_e$ is an optimal solution of the problem in (\ref{problem1})

Let $z_1 \in\mathbb{R}^n$ be any constant vector and $\mathbf{z}_1=1_N\otimes z_1$. 

Define a function $V_{\mathbf{z}_1}\in\mathbb{R}$ as
\begin{eqnarray}
V_{\mathbf{z}_1}(t,\mathbf{x},\dot{\mathbf{x}})&=&\frac{1}{2}(\mathbf{x}+e^{{-\alpha}_t} \dot{\mathbf{x}}-\mathbf{z}_1)^T(\mathbf{x}+e^{{-\alpha}_t} \dot{\mathbf{x}}-\mathbf{z}_1)\notag \\
& &+\frac{1}{2}{(\mathbf{x}-\mathbf{z}_1)}^T{(\mathbf{x}-\mathbf{z}_1)}+\frac{k_1}{2}\mathbf{x}^T\tilde{L}\mathbf{x} \notag \\
& &+e^{\beta_t}(F_t(\mathbf{x})-F_t(\mathbf{z}_1)).
\end{eqnarray}

Then,
\begin{eqnarray}
V_{\mathbf{z}_1}(t,\mathbf{x},\dot{\mathbf{\mathbf{x}}})&=&\frac{1}{2}(\mathbf{x}+e^{{-\alpha}_t} \dot{\mathbf{x}}-\mathbf{z}_1)^T(\mathbf{x}+e^{{-\alpha}_t} \dot{\mathbf{x}}-\mathbf{z}_1)\notag \\
& &+\frac{1}{2}{(\mathbf{x}-\mathbf{z}_1)}^T{(\mathbf{x}-\mathbf{z}_1)}+\frac{k_1}{2}\mathbf{x}^T\tilde{L}\mathbf{x}  \notag \\
& &+e^{\beta_t}(F_t(\mathbf{x})-F_t(\mathbf{z}_0(t))) \notag \\
& &+e^{\beta_t}(F_t(\mathbf{z}_0(t))-F_t(\mathbf{z}_1)) \notag \\
&\geq&\frac{1}{2}{(\mathbf{x}-\mathbf{z}_1)}^T{(\mathbf{x}-\mathbf{z}_1)} \notag \\
& &+e^{\beta_t}(F_t(\mathbf{z}_0(t))-F_t(\mathbf{z}_1)) \notag \\
&\geq& e^{\beta_t}(F_t(\mathbf{z}_0(t))-F_t(\mathbf{z}_1)), \label{Lya_2}
\end{eqnarray}
where $\mathbf{z}_0(t)=[\mathbf{x}_{1}^{\dagger T}(t),\cdots,\mathbf{x}_{N}^{\dagger T}(t)]^T$.

Taking the derivative of $V_{\mathbf{z}_1}(t,\mathbf{x},\dot{\mathbf{x}})$ gives
\begin{eqnarray}
\dot{V}_{\mathbf{z}_1}&=&(\mathbf{x}+e^{{-\alpha}_t} \dot{\mathbf{x}}-\mathbf{z}_1)^T(\dot{\mathbf{x}}-\dot{\alpha}_te^{{-\alpha}_t}\dot{\mathbf{x}}+e^{{-\alpha}_t} \ddot{\mathbf{x}}) \notag \\
& &+{(\mathbf{x}-\mathbf{z}_1)}^T\dot{\mathbf{x}}+k_1\mathbf{x}^T\tilde{L}\dot{\mathbf{x}} \notag \\
& &+\dot{\beta}_te^{\beta_t}(F_t(\mathbf{x})-F_t(\mathbf{z}_1))
\notag \\
& &+e^{\beta_t}(\dot{\mathbf{x}}^T\nabla F_t(\mathbf{x})+\nabla_t F_t(\mathbf{x})-\nabla_t F_t(\mathbf{z_1})) \notag \\
&=&(\mathbf{x}+e^{{-\alpha}_t} \dot{\mathbf{x}}-\mathbf{z}_1)^T(-e^{\alpha_t+\beta_t}\nabla F_t(\mathbf{x}) \notag \\
& &-\dot{\mathbf{x}}-k_1e^{{\alpha}_t}\tilde{L}\mathbf{x})+{(\mathbf{x}-\mathbf{z}_1)}^T\dot{\mathbf{x}}+k_1\mathbf{x}^T\tilde{L}\dot{\mathbf{x}} \notag \\
& &+\dot{\beta}_te^{\beta_t}(F_t(\mathbf{x})-F_t(\mathbf{z}_1))+e^{\beta_t}\dot{\mathbf{x}}^T\nabla F_t(\mathbf{x})  \notag \\
& &+e^{\beta_t} (\nabla_t F_t(\mathbf{x})-\nabla_t F_t(\mathbf{z_1}))  \notag \\
&\leq&-(e^{\alpha_t+\beta_t}-\dot{\beta}_te^{\beta_t})(F_t(\mathbf{x})-F_t(\mathbf{z}_1))  \notag \\
& &+e^{\beta_t}(\nabla_t F_t(\mathbf{x})-\nabla_t F_t(\mathbf{z_1})) \notag \\
& &-e^{-\alpha_t}\dot{\mathbf{x}}^T\dot{\mathbf{x}}-k_1e^{\alpha_t}\mathbf{x}^T\tilde{L}\mathbf{x},   \label{derivative2}
\end{eqnarray} 
where in the last step we use the fact that $\tilde{L}\mathbf{z}_1=0$.

Similarly, let (\ref{f}), (\ref{eq}), and (\ref{fff}) hold. Letting $z_1=\tilde{x}(T)$ and $\mathbf{\tilde{x}}(T)=1_{N}\otimes\tilde{x}(T)$, it holds that  
\begin{eqnarray}
& & \int_{T_0}^T(\sigma(F_t(\mathbf{x}(t))-F_t(\mathbf{\tilde{x}}(T)))+k_1 e^{\alpha_t}\mathbf{x}^T(t)\tilde{L}\mathbf{x}(t)) dt  \notag \\
&\leq&V_{\mathbf{\tilde{x}}}(T_0,\mathbf{x}(T_0),\dot{\mathbf{x}}(T_0)) \notag \\
& &	-\frac{1}{2}{(\mathbf{x}(T)-\mathbf{\tilde{x}}(T))}^T{(\mathbf{x}(T)-\mathbf{\tilde{x}}(T))} \notag \\
& &-e^{\beta_T}(F_t(\mathbf{z}_0(T))-F_t(\mathbf{\tilde{x}}(T))) \notag \\
& &+\int_{T_0}^{T}e^{\beta_t}(\nabla_t F_t(\mathbf{x}(t))-\nabla_t F_t(\mathbf{\tilde{x}}(T)))dt  \label{equation2aa} \\
&\leq& V_{\mathbf{\tilde{x}}}(T_0,\mathbf{x}(T_0),\dot{\mathbf{x}}(T_0))-e^{\beta_T}(F_t(\mathbf{z}_0(T))-F_t(\mathbf{\tilde{x}}(T))) \notag \\
& &+\int_{T_0}^{T}e^{\beta_t}(\nabla_t F_t(\mathbf{x}(t))-\nabla_t F_t(\mathbf{\tilde{x}}(T)))dt. \label{equation2a}
\end{eqnarray} 

Based on Assumption \ref{a2dis}, \ref{as} and \ref{a32}, there exists a constant $\varrho>0$ independent of $T$ such that 
\begin{eqnarray}
& &F_t(\mathbf{x}(t))-F_t(\mathbf{\tilde{x}}(T)) \notag \\
&\geq& \text{min}_{\mathbf{x}\in\mathbb{R}^{Nn}}F_t(\mathbf{x})-F_t(\mathbf{\tilde{x}}(T))\geq-\varrho. 
\end{eqnarray} 

Then, 
\begin{eqnarray}
& & \int_{T_0}^T(-\sigma\varrho+k_1e^{\alpha_t}\mathbf{x}^T(t)\tilde{L}\mathbf{x}(t)) dt  \notag \\
&\leq&V_{\mathbf{\tilde{x}}}(T_0,\mathbf{x}(T_0),\dot{\mathbf{x}}(T_0))-e^{\beta_T}(F_t(\mathbf{z}_0(T))-F_t(\mathbf{\tilde{x}}(T))) \notag \\
& &+\int_{T_0}^{T}e^{\beta_t}(\nabla_t F_t(\mathbf{x}(t))-\nabla_t F_t(\mathbf{\tilde{x}}(T)))dt. \label{eq2}
\end{eqnarray} 

Similarly to Section \ref{S3_1}, according to Assumptions \ref{a1dis}-\ref{boundedness1}, it can be verified that if $m< -1$,
%\begin{eqnarray}
%\int_{T_0}^T(w(t)\mathbf{x}^T\tilde{L}\mathbf{x}) dt\leq \frac{\sigma\varrho}{k_1}(T-T_0)+ \frac{\varpi}{ k_1}  
%\end{eqnarray}
\begin{eqnarray}
\int_{T_0}^T(e^{\alpha_t}\mathbf{x}^T(t)\tilde{L}\mathbf{x}(t)-\frac{\sigma\varrho}{k_1}) dt\leq \frac{\varpi}{ k_1}   \label{c_1}
\end{eqnarray}
with $ \varpi$ being some positive constant independent of $T$.

Suppose that
\vspace{-4pt}
\begin{eqnarray}
b_0\geq 1. \label{pb}
\end{eqnarray}

Then, 
\begin{eqnarray}
\int_{T_0}^Te^{\alpha_t}\mathbf{x}^T(t)\tilde{L}\mathbf{x}(t) dt\leq \frac{\sigma\varrho}{k_1}(T-T_0)+\frac{\varpi}{ k_1}\triangleq \epsilon_1(T),   \label{c_2}
\end{eqnarray}
and
\begin{eqnarray}
\int_{T_0}^T\mathbf{x}^T(t)\tilde{L}\mathbf{x}(t) dt&\leq& \frac{\varrho}{k_1b_0^{-m}}(T-T_0)+\frac{\varpi}{ k_1\sigma b_0^{-m}} \notag \\
&\triangleq& \epsilon_2(T),   \label{c_3}
\end{eqnarray}
where we use the fact that $e^{\alpha_t}\geq\sigma b_0 ^{-m}$.

Then, for any $i,j\in\mathcal{V}$,
\begin{eqnarray}
\int_{T_0}^Te^{\alpha_t}\left\Vert\mathbf{x}_i(t)-\mathbf{x}_j(t)\right\Vert^2 dt &\leq&
N\int_{T_0}^Te^{\alpha_t}\mathbf{x}^T(t)\tilde{L}\mathbf{x}(t) dt   \notag \\
&\leq& \frac{N}{w_\text{min}} \epsilon_1(T),
\end{eqnarray}
where $w_\text{min}>0$ is the minimal weight of the graph. 

According to the convexity of the square function, it holds that
\begin{eqnarray}
& &(\int_{T_0}^Te^{\frac{\alpha_t}{2}}\left\Vert\mathbf{x}_i(t)-\mathbf{x}_j(t)\right\Vert dt)^2 \notag \\
&\leq& \int_{T_0}^Te^{\alpha_t}\left\Vert\mathbf{x}_i(t)-\mathbf{x}_j(t)\right\Vert^2 dt 
\leq \frac{N}{w_\text{min}}\epsilon_1(T).
\end{eqnarray}

Then, 
\begin{eqnarray}
\int_{T_0}^Te^{\frac{\alpha_t}{2}}\left\Vert\mathbf{x}_i(t)-\mathbf{x}_j(t)\right\Vert dt
&\leq& \sqrt{ \frac{N}{w_\text{min}}\epsilon_1(T)}.
\end{eqnarray}

Similarly, the following inequality holds
\begin{eqnarray}
\int_{T_0}^T\left\Vert\mathbf{x}_i(t)-\mathbf{x}_j(t)\right\Vert dt
&\leq& \sqrt{ \frac{N}{w_\text{min}}\epsilon_2(T)}.
\end{eqnarray}

According to (\ref{equation2aa}), it can be obtained that 
\begin{eqnarray}
\left\Vert\mathbf{x}(T)-\mathbf{\tilde{x}}(T)\right\Vert^2 \leq 2 (\varpi+\sigma\varrho(T-T_0)).
\end{eqnarray}

It follows that
\begin{eqnarray}
\left\Vert\mathbf{x}(T)\right\Vert \leq \left\Vert\mathbf{\tilde{x}}(T)\right\Vert+\sqrt{2 \left\vert\varpi-\sigma\varrho T_0\right\vert}+\sqrt{2\sigma\varrho T}. \label{aa2}
\end{eqnarray}

Since (\ref{aa2}) hold for any $T>T_0$, and $\left\Vert\mathbf{\tilde{x}}(T)\right\Vert$ is upper bounded and the bound is independent of $T$, it holds that for any $T_0<t\leq T$, 
\begin{eqnarray}
\left\Vert\mathbf{x}(t)\right\Vert \leq \tilde{c}+\sqrt{2 \left\vert\varpi-\sigma\varrho T_0\right\vert}+\sqrt{2\sigma\varrho t}, \label{aa1}
\end{eqnarray}
where $\tilde{c}$ is the bound of $\left\Vert\mathbf{\tilde{x}}(T)\right\Vert$, which is independent of $ T$  .

According to Assumption \ref{boundednessgra}, there exist $a_2$ and $b_2$ such that
\begin{eqnarray}
\left\Vert \nabla f_{i,t}(\mathbf{x}_i(t)) \right\Vert \leq a_2 t^{m_0}+b_2.
\end{eqnarray}

Suppose that
\begin{eqnarray}
m\leq -2m_0. \label{pm}
\end{eqnarray}
It can be obtained that for any $t>T_0$,
\begin{eqnarray}
e^{\frac{\alpha_t}{2}}\geq \sqrt{\sigma}t^{m_0}.
\end{eqnarray}

Define 
\begin{eqnarray}
\mathcal{R}_{si}(T)=\int_{T_0}^T (\sum_{j=1}^{N}f_{j,t}(\mathbf{x}_i(t))-\sum_{j=1}^{N}f_{j,t}(\tilde{x}(T)))dt.
\end{eqnarray}

Then,
\begin{eqnarray}
\mathcal{R}_{si}(T)
&=&\int_{T_0}^T(\sum_{j=1}^{N}f_{j,t}(\mathbf{x}_j(t))-\sum_{j=1}^{N}f_{j,t}(\tilde{x}(T)))dt \notag \\
& &-\int_{T_0}^T (\sum_{j=1}^{N}f_{j,t}(\mathbf{x}_j(t))-\sum_{j=1}^{N}f_{j,t}(\mathbf{x}_i(t)))dt \notag \\
&\leq&\int_{T_0}^T(F_t(\mathbf{x}(t))-F_t(\mathbf{\tilde{x}}(T)))dt \notag \\
& &+\int_{T_0}^T \sum_{j=1}^{N} \left\Vert \nabla f_{j,t}(\mathbf{x}_j(t)) \right\Vert \left\Vert\mathbf{x}_i(t)-\mathbf{x}_j(t)\right\Vert dt \notag \\
&\leq& \frac{\varpi}{\sigma}+\int_{T_0}^T \sum_{j=1}^{N} (a_2 t^{m_0}+b_2) \left\Vert\mathbf{x}_i(t)-\mathbf{x}_j(t)\right\Vert dt \notag \\
&\leq& \frac{\varpi}{\sigma}+ \int_{T_0}^T  \sum_{j=1}^{N} (\frac{a_2}{\sqrt{\sigma}}e^{\frac{\alpha_t}{2}}+b_2)\left\Vert\mathbf{x}_i(t)-\mathbf{x}_j(t)\right\Vert dt \notag \\
&\leq&\frac{\varpi}{\sigma}+\frac{Na_2}{\sqrt{\sigma}}\sqrt{ \frac{N}{w_\text{min}}\epsilon_1(T)}+Nb_2\sqrt{ \frac{N}{w_\text{min}}\epsilon_2(T)}, \notag 
\end{eqnarray}
which implies that
\begin{eqnarray}
\mathcal{R}_{s,dis}(T)\leq\frac{\varpi}{\sigma}+\frac{Na_2}{\sqrt{\sigma}}\sqrt{ \frac{N}{w_\text{min}}\epsilon_1(T)}+Nb_2\sqrt{ \frac{N}{w_\text{min}}\epsilon_2(T)}. \notag 
\end{eqnarray}

\begin{theorem}
	Assume that Assumptions \ref{a1dis}-\ref{boundednessgra}  hold. The DOCT-NAG algorithm in (\ref{algorithm2}) satisfies that for all $T>T_0$,
	\begin{eqnarray}
	\mathcal{R}_{s,dis}(T)&\leq& \frac{\varpi}{\sigma}+\frac{Na_2}{\sqrt{\sigma}}\sqrt{ \frac{N}{w_\text{min}}\epsilon_1(T)} \notag \\ & &+Nb_2\sqrt{ \frac{N}{w_\text{min}}\epsilon_2(T)}, 
   \label{theore}
	\end{eqnarray}
	provided that the conditions in  (\ref{f}), (\ref{eq}), (\ref{fff}), (\ref{pb}) and (\ref{pm}) hold and $m<-1$, where $\sigma$ is any positive parameter, $w_\text{min}>0$ is the minimal weight, $a_2$ and $a_3$ are defined in Assumption \ref{boundednessgra}, $\varpi$ is defined in (\ref{c_1}), and $\epsilon_1(T)$ and $\epsilon_2(T)$ are defined in (\ref{c_2}) and (\ref{c_3}). 
	 \label{theorem4}
\end{theorem}

\begin{remark}
The proposed distributed algorithm in (\ref{algorithm2}) uses a proportional consensus term rather than a proportional-integral term, which benefits from the introduction of time-varying gains\cite{liu2014continuous}. The time-varying parameters behave like a penalty parameter which goes to infinity with time.
\end{remark}

\subsection{Dynamic Regret Analysis}

Similar to Section \ref{aaa}, we have the following conclusion.

\begin{theorem}
Assume that Assumptions \ref{a1dis}-\ref{boundednessgra}  hold. The DOCT-NAG algorithm in (\ref{algorithm2}) satisfies that for all $T>T_0$,
	\begin{eqnarray}
	\mathcal{R}_{d,dis}(T)&\leq& c_4(T-T_0)+\frac{\varpi}{\sigma}+\frac{Na_2}{\sqrt{\sigma}}\sqrt{ N\epsilon_1(T)} \notag \\
	& &+Nb_2\sqrt{ N\epsilon_2(T)}, \label{theoremconditionab}
	\end{eqnarray}
	provided that the conditions in  (\ref{f}), (\ref{eq}), (\ref{fff}), (\ref{pb}) and (\ref{pm}) hold and $m<-1$, where $\sigma$ is any positive parameter, $w_\text{min}>0$ is the minimal weight, $a_2$ and $a_3$ are defined in Assumption \ref{boundednessgra}, $\varpi$ is defined in (\ref{c_1}), and $\epsilon_1(T)$ and $\epsilon_2(T)$ are defined in (\ref{c_2}) and (\ref{c_3}), and $c_4$ is a positive constant independent of $T$.
\end{theorem}

\section{Simulation} \label{S4} 

\subsection{OCT-NAG} \label{simu1}

Consider the following online convex optimization problem:
\begin{eqnarray}
\text{min}_{x\in\mathbb{R}} f_t(x)=x^2+\text{sin}(t)\text{sin}(x). \label{problemsimu}
\end{eqnarray}

The function is strictly convex at any time since its second order derivative is positive definite. Furthermore,
\begin{eqnarray}
\nabla f_t(x_t^*)=2x^*_t+\text{sin}(t)\text{cos}(x^*_t)=0,
\end{eqnarray}
and $x^*_t$ exists and is uniformly bounded. In addition, for any $T>T_0$,
\begin{eqnarray}
\int_{T_0}^T f_t(x) dt&=&\int_{T_0}^T (x^2+\text{sin}(t)\text{sin}(x))dt \notag \\
&=&x^2(T-T_0)-\text{sin}(x)\text{cos}(t) |_{T_0}^T,
\end{eqnarray}
and 
\begin{eqnarray}
\nabla\int_{T_0}^T f_t(x) dt&=&2(T-T_0)x-(\text{cos}(T) \notag \\ & &-\text{cos}(T_0))\text{cos}(x), \label{me}
\end{eqnarray}
and 
\begin{eqnarray}
\nabla^2\int_{T_0}^T f_t(x) dt&=&2(T-T_0)+(\text{cos}(T) \notag \\ & &-\text{cos}(T_0))\text{sin}(x).
\end{eqnarray}

Since the function $g(\gamma)=2\gamma+\text{cos}(\gamma)\text{sin}(x)$ is a strictly monotonically increasing function on $\gamma$ for any $x$, we can obtain that $\nabla^2\int_{T_0}^T f_t(x) dt>0$, which implies that $\int_{T_0}^T f_t(x) dt$ is strictly convex. Moreover, from (\ref{me}), it can be verified that the static optimal solution $\tilde{x}(T)$ satisfying $\nabla\int_{T_0}^T f_t(x) dt=0$ exists and satisfies 
\begin{eqnarray}
\left\vert \tilde{x}(T) \right\vert=\left\vert\frac{\text{cos}(T)-\text{cos}(T_0)}{2(T-T_0)}\text{cos}(\tilde{x}(T))\right\vert \notag \\
=\left\vert-\frac{\text{sin}(\frac{T+T_0}{2})\text{sin}(\frac{T-T_0}{2})}{T-T_0}\text{cos}(\tilde{x}(T))\right\vert< \frac{1}{2}, \label{me1}
\end{eqnarray}
where we use the fact that $\left\vert\frac{\text{sin}(a)}{a}\right\vert<1$ for $a>0$, which implies that $\tilde{x}(T)$ is bounded with the bound independent of $T$. Then, it can be verified that Assumptions \ref{a1}-\ref{boundedness} hold.

Matlab Simulink is used in the simulation and the solver is chosen as ode23tb. For better demonstration of both the algorithm evolution and the final algorithm regret and errors, we select the optimal solution $x^*_t |_{t=0}=0$ at $t=0$ as the initial point. For other initial points, the conclusions and analysis in this section still hold. In addition, it is found that a large $|m|$ may lead to high stiffness of the system, which may affect the performance of the algorithm in the simulation.

\begin{figure}[htb]
	\centering
	\includegraphics[width=7cm]{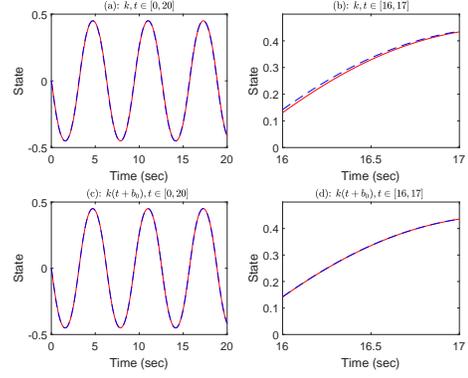} \caption{The states when $e^{\alpha_t+\beta_t}-\dot{\beta}_te^{\beta_t}=\sigma$ and $e^{\alpha_t+\beta_t}-\dot{\beta}_te^{\beta_t}=\sigma(t+b_0)$. The dashed blue line represents $x^*_t$ and the red line represents the state.}
	\label{ch22state}
\end{figure}

Consider the algorithm in (\ref{algorithm}) and let $m=-20, \sigma=20$, $b_0=2$. Fig. 	\ref{ch22state} (a), (b) shows the comparison of the state variable with the time-varying optimal solution $x^*_t$ when $e^{\alpha_t+\beta_t}-\dot{\beta}_te^{\beta_t}=\sigma$ and Fig. \ref{ch22state}  (c), (d) shows the result when $e^{\alpha_t+\beta_t}-\dot{\beta}_te^{\beta_t}=\sigma(t+b_0)$. Both the algorithms track the optimal solution, but it can be seen that the tracking for the parameter satisfying  $e^{\alpha_t+\beta_t}-\dot{\beta}_te^{\beta_t}=\sigma(t+b_0)$ is more precise.

\subsubsection{Static Regret}

\begin{figure}[htb]
	\centering
	\includegraphics[width=7cm]{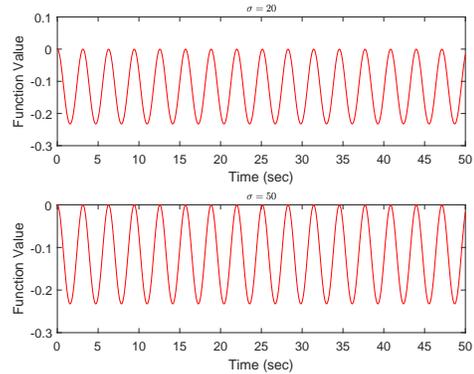} \caption{Function values when $\sigma=20$ and $\sigma=50$.}
	\label{ch21}
\end{figure}

\begin{figure}[htb]
	\centering
	\includegraphics[width=7cm]{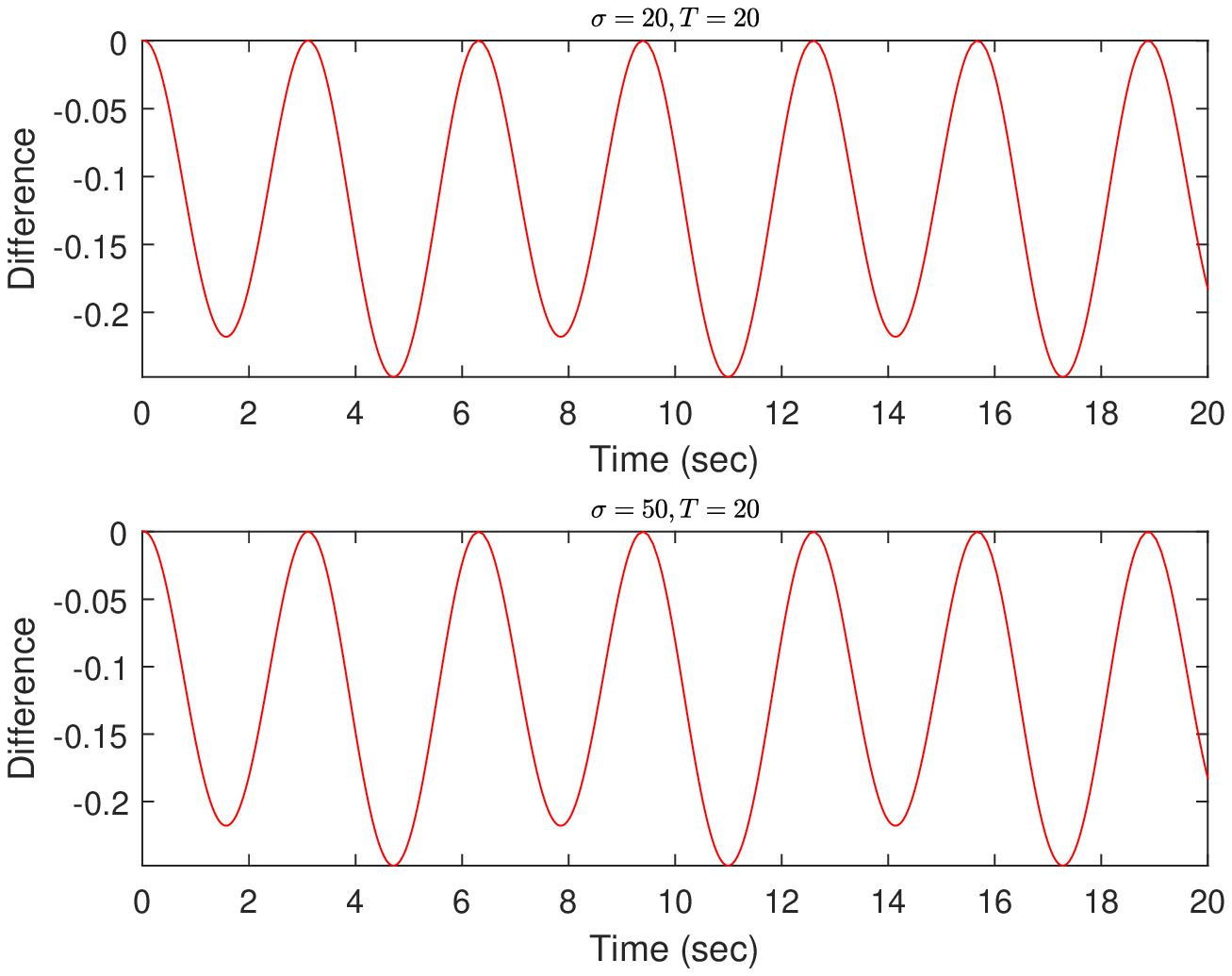} \caption{The differences with the static optimal function value when $\sigma=20, T=20$ and $\sigma=50, T=20$.}
	\label{ch22}
\end{figure}

\begin{figure}[htb]
	\centering
	\includegraphics[width=7cm]{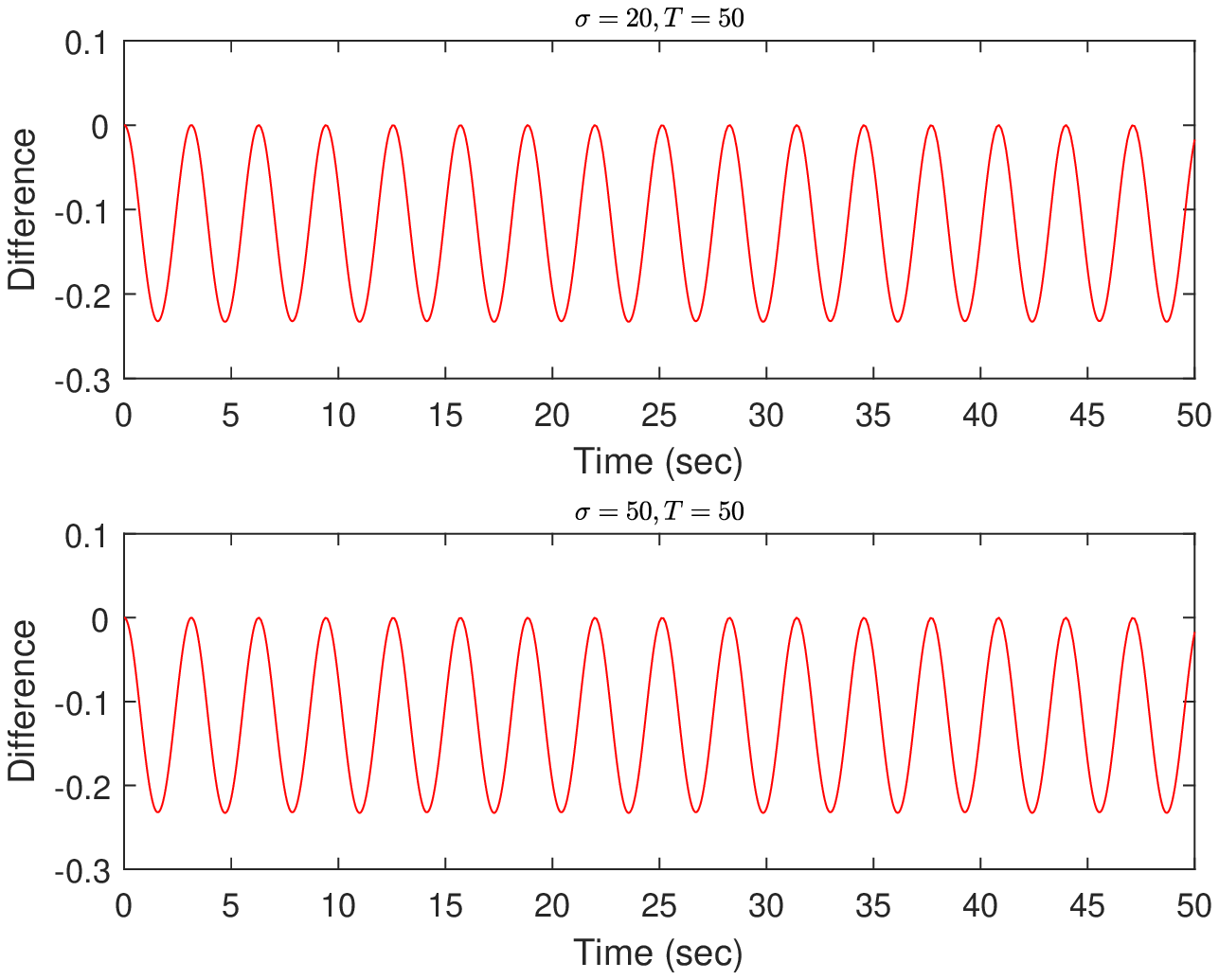} \caption{The differences with the static optimal function value when $\sigma=20, T=50$ and $\sigma=50, T=50$.}
	\label{ch24}
\end{figure}

Let (\ref{algorithm}), (\ref{f}), (\ref{eq}), and (\ref{fff}) be the updating law with $m=-20, \sigma=20$, $b_0=2$, and $m=-50, \sigma=50$, $b_0=2$, respectively. The simulation results are shown in Figs. \ref{ch21}-\ref{ch24}. Fig. \ref{ch21} shows the function values calculated by the algorithm when $\sigma=20$ and $\sigma=50$. Fig. \ref{ch22} shows the differences of the algorithm function values with the static optimal function value when $\sigma=20$ and $\sigma=50$ for $T=20$, i.e., $f_t(x(t))-f_t(\tilde{x}(20))$. Fig. \ref{ch24} shows the differences of the algorithm function values with the static optimal function value when $\sigma=20$ and $\sigma=50$ for $T=50$, i.e., $f_t(x(t))-f_t(\tilde{x}(50))$. The curves are always in the negative half plane. Thus, the conclusion in Theorem \ref{theorem1} is verified. 

\begin{figure}[htb]
	\centering
	\includegraphics[width=7cm]{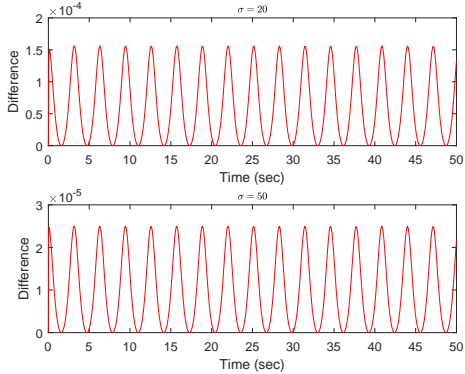} \caption{The differences with the dynamic optimal function value when $\sigma=20$ and $\sigma=50$.}
	\label{ch23}
\end{figure}

\begin{figure}[htb]
	\centering
	\includegraphics[width=7cm]{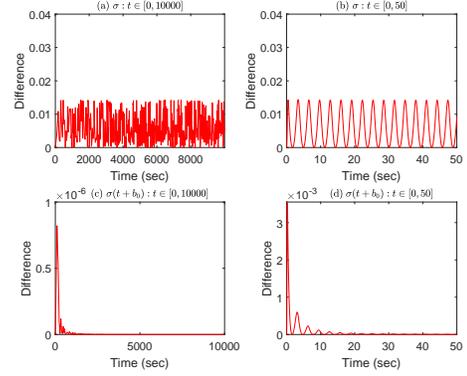} \caption{The differences with the dynamic optimal function value for $e^{\alpha_t+\beta_t}-\dot{\beta}_te^{\beta_t}=\sigma$ and $e^{\alpha_t+\beta_t}-\dot{\beta}_te^{\beta_t}=\sigma(t+b_0)$.}
	\label{ch25}
\end{figure}

\begin{figure}[htb]
	\centering
	\includegraphics[width=7cm]{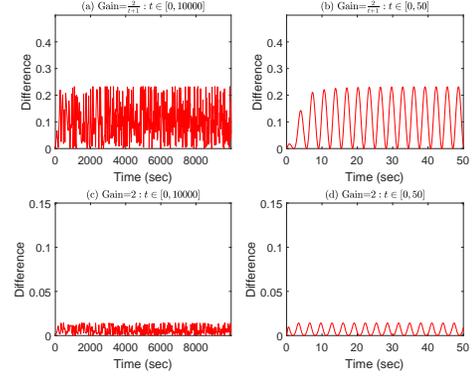} \caption{The differences with the dynamic optimal function value for the continuous-time Online Gradient Descent Methods with a time-varying gain and a time-invariant gain.}
	\label{ch26}
\end{figure}

\begin{figure}[htb]
	\centering
	\includegraphics[width=7cm]{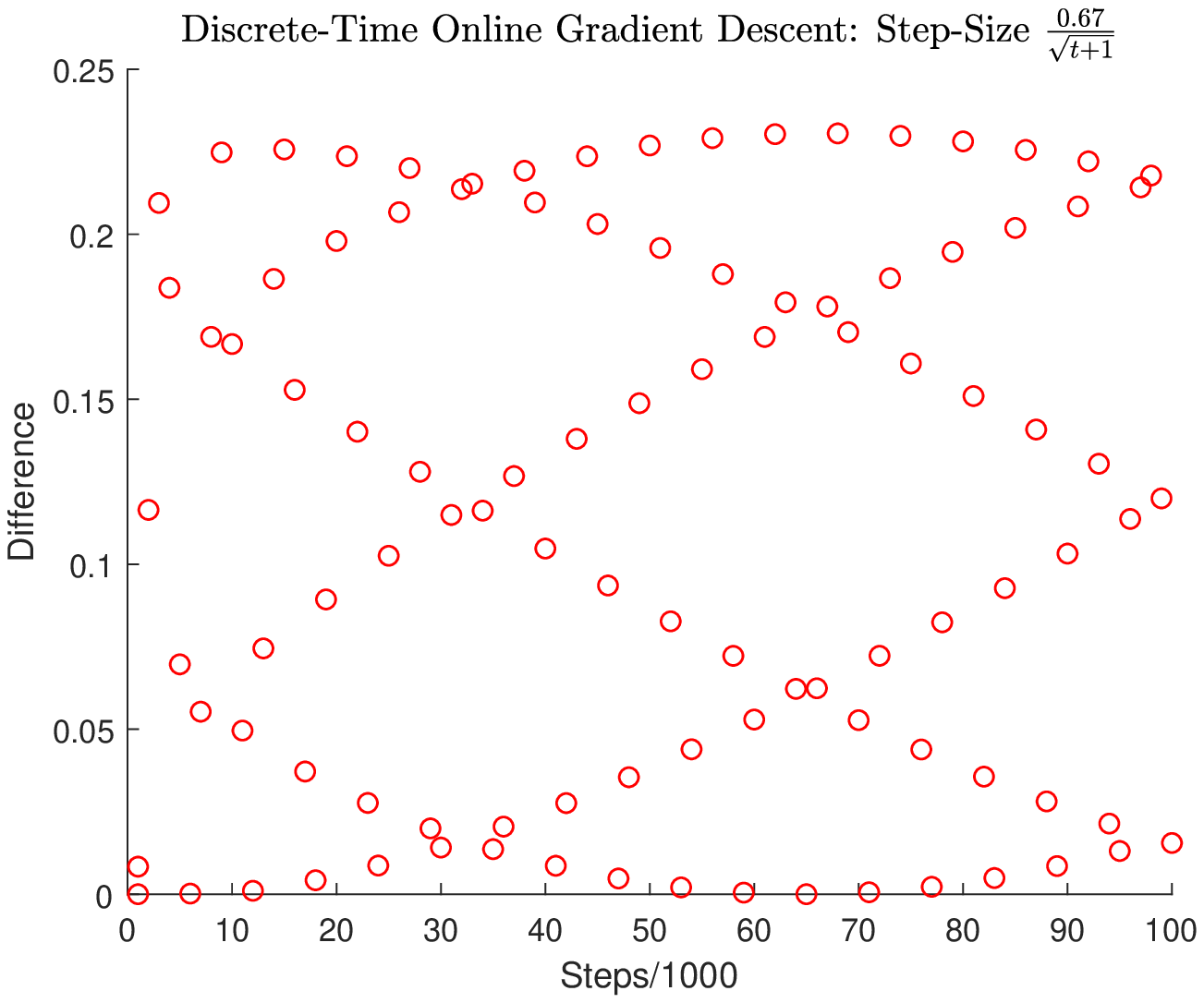} \caption{The difference with the dynamic optimal function value for the discrete-time Online Gradient Descent Method with a step-size $\frac{0.67}{\sqrt{t+1}}$.}
	\label{ch27}
\end{figure}

\begin{figure}[htb]
	\centering
	\includegraphics[width=7cm]{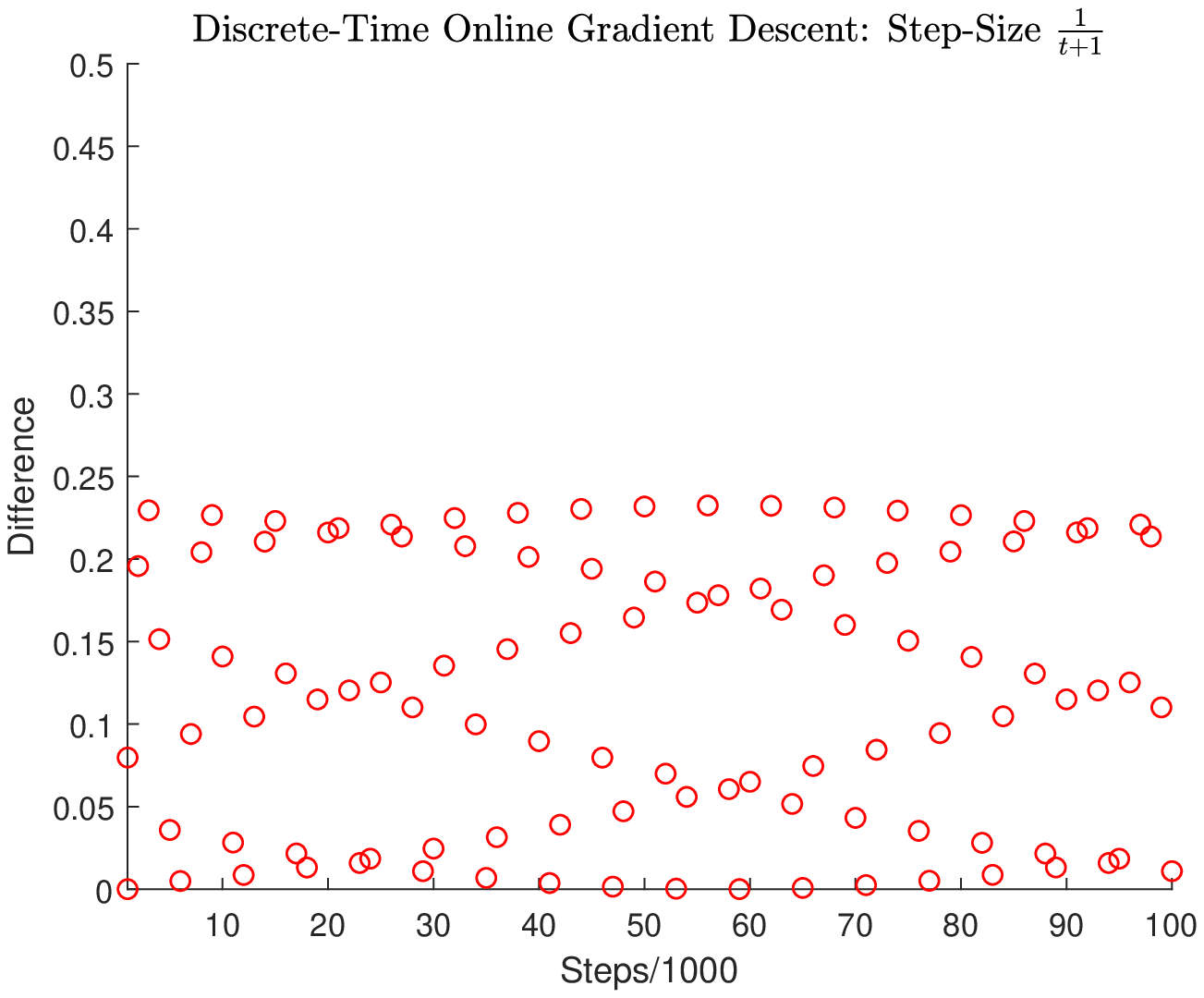} \caption{The difference with the dynamic optimal function value for the discrete-time Online Gradient Descent Method with a step-size $\frac{1}{{t+1}}$.}
	\label{ch28}
\end{figure}

\begin{figure}[htb]
	\centering
	\includegraphics[width=7cm]{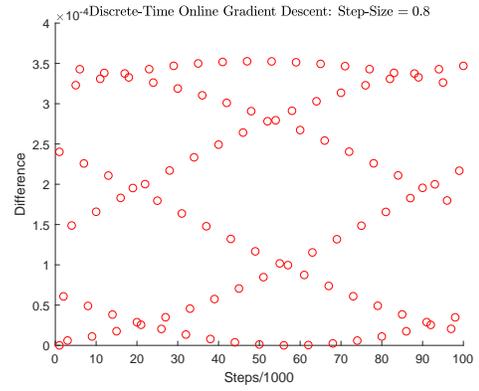} \caption{The difference with the dynamic optimal function value for the discrete-time Online Gradient Descent Method with a fixed step-size $0.8$.}
	\label{ch29}
\end{figure}

\begin{figure}[htb]
	\centering
	\includegraphics[width=7cm]{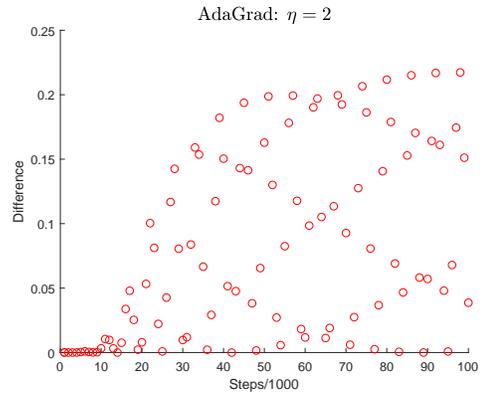} \caption{The difference with the dynamic optimal function value for AdaGrad with $\eta=2$.}
	\label{ch210}
\end{figure}

\begin{figure}[htb]
	\centering
	\includegraphics[width=7cm]{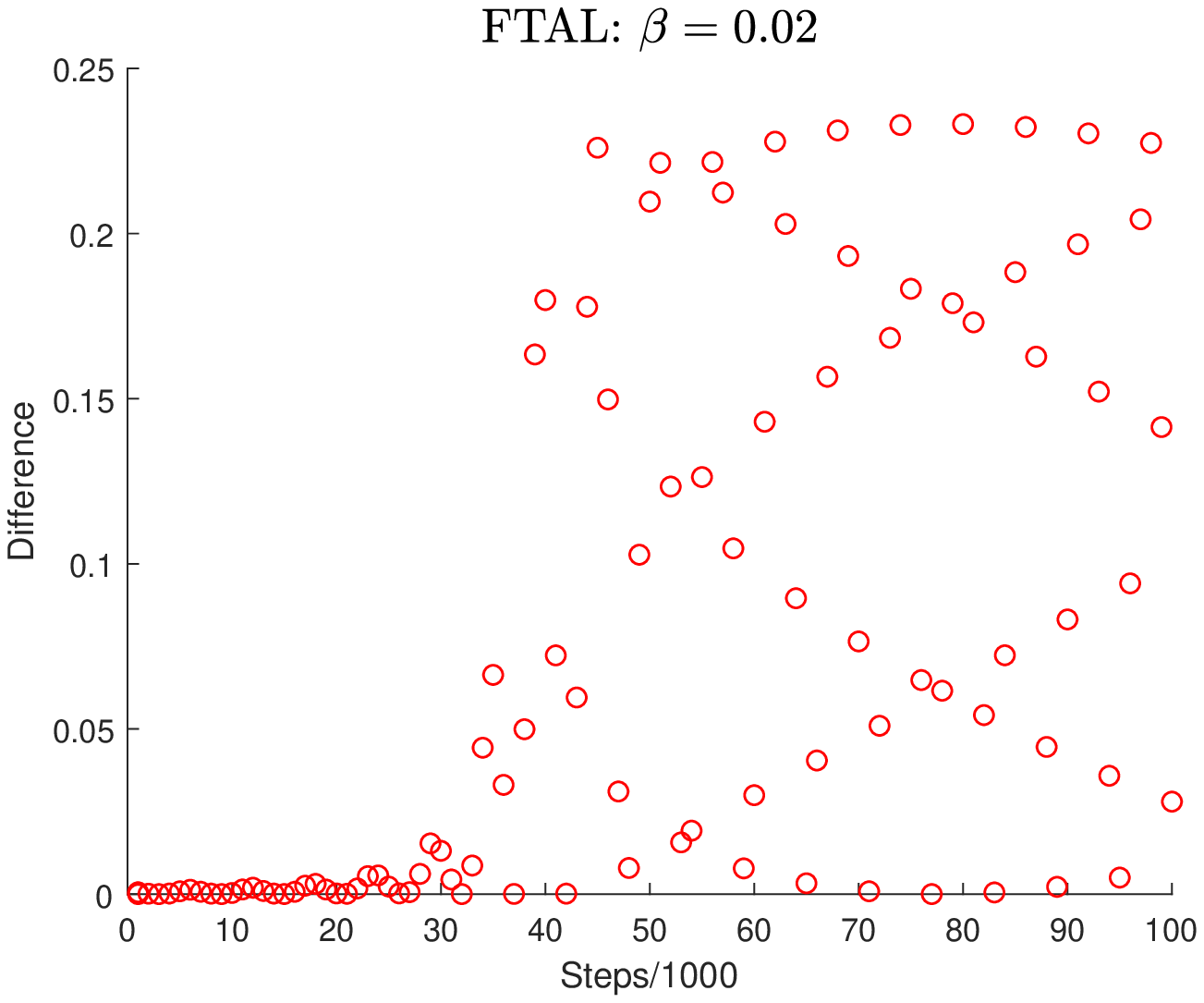} \caption{The difference with the dynamic optimal function value for FTAL with $\beta=0.02$.}
	\label{ch211}
\end{figure}

\subsubsection{Dynamic Regret}

We compare the real-time performance of the algorithm with other algorithms according to the difference of the algorithm function value with the dynamic optimal function value i.e., $f_t(x(t))-f_t(x^*_t)$. Since it is difficult to obtain an explicit expression of the time-varying optimal solution of the function, in the following comparisons, we select enough discrete-time points and sample the algorithms at these points for comparisons.

Letting (\ref{algorithm}), (\ref{f}), (\ref{eq}), and (\ref{fff}) be the updating law, Fig. \ref{ch23} shows the differences of the algorithm function values with the dynamic optimal function value when $\sigma=20$ and $\sigma=50$. It can be seen that the error is smaller when $\sigma$ is larger.

In the following simulation, we select similar parameters for the algorithms for more precise comparisons. The bound of the gradients is selected to be $G=3$ and the diameter of the convex set is selected to be $D=2$ \cite{hazan2016introduction, hazan2007logarithmic}. 

For the algorithm in  (\ref{algorithm}), (\ref{f}), (\ref{eq}), and (\ref{fff}) which corresponds to the case where $e^{\alpha_t+\beta_t}-\dot{\beta}_te^{\beta_t}=\sigma$, Figs. \ref{ch25}(a) and \ref{ch25}(b) show the algorithm results for $10000s$ and $50s$, respectively. For the algorithm in (\ref{algorithm}), (\ref{eq}), (\ref{f_2}),  and (\ref{fff_2}), which corresponds to the case where $e^{\alpha_t+\beta_t}-\dot{\beta}_te^{\beta_t}=\sigma(t+b_0)$, Figs. \ref{ch25}(c) and \ref{ch25}(d) show the results. The parameters are selected to be $m=-2, \sigma=2, b_0=2$. It can be seen that the algorithm with $e^{\alpha_t+\beta_t}-\dot{\beta}_te^{\beta_t}=\sigma(t+b_0)$ has a better performance.

Fig. \ref{ch26} shows the differences by using the continuous-time Online Gradient Descent method, where Figs. \ref{ch26}(a) and \ref{ch26}(b) show the case where the gain is selected to be $\frac{2}{t+1}$ and Figs. \ref{ch26}(c) and \ref{ch26}(d) show the case where the gain is fixed to be $2$. 

Fig. \ref{ch27} shows the difference by using the discrete-time Online Gradient Descent method \cite{hazan2016introduction} with a step-size $\frac{0.67}{\sqrt{t+1}}$ ($t=0, 1, \cdots$). The sampling period for the algorithm is $0.1$s. The algorithm runs for 100000 steps and a sampling point is selected in every 1000 points to be compared with the dynamic optimal function value.

Fig. \ref{ch28} shows the difference by using the discrete-time Online Gradient Descent method \cite{hazan2016introduction} with a step-size $\frac{1}{t+1}.$ 

Fig. \ref{ch29} shows the difference by using the discrete-time Online Gradient Descent method \cite{hazan2016introduction} with a fixed step-size $0.8$. According to \cite{hazan2016introduction}, this case can be viewed as a special case of the Online Mirror Descent algorithm in an agile version. Thus, this figure can be viewed as a comparison with the Online Mirror Descent method.

Fig. \ref{ch210} shows the difference by using the AdaGrad method \cite{hazan2016introduction} with parameter $\eta=2$.

Fig. \ref{ch211} shows the difference by using the FTAL method \cite{hazan2007logarithmic} with parameter $\beta=0.02$.

It can be seen that for the proposed method with a scaling condition $e^{\alpha_t+\beta_t}-\dot{\beta}_te^{\beta_t}=\sigma(t+b_0)$, the upper bound of the difference with the dynamic optimal function value is persistently decreasing (approximately exponentially or polynomially). This phenomenon doesn't occur in other algorithms. The simulation verifies the effectiveness and efficiency of the proposed algorithm. 

\subsection{DOCT-NAG}

In this section, we consider a group of 6 agents that coordinate with each other to solve an optimization problem
\begin{eqnarray}
\text{min } f_t(x)=\sum_{i=1}^6 f_{i,t}(x), \label{systemsim}
\end{eqnarray}
where the decision variable $x=[x_1, \cdots, x_6]^T\in\mathbb{R}^6$, and 
\begin{eqnarray}
f_{i,t}(x)= 10(x^Tx+0.2\times i \times \text{sin}(t)\text{sin}(x_i)). 
\end{eqnarray}

\begin{figure}[htb]
	\centering
	\includegraphics[width=6cm]{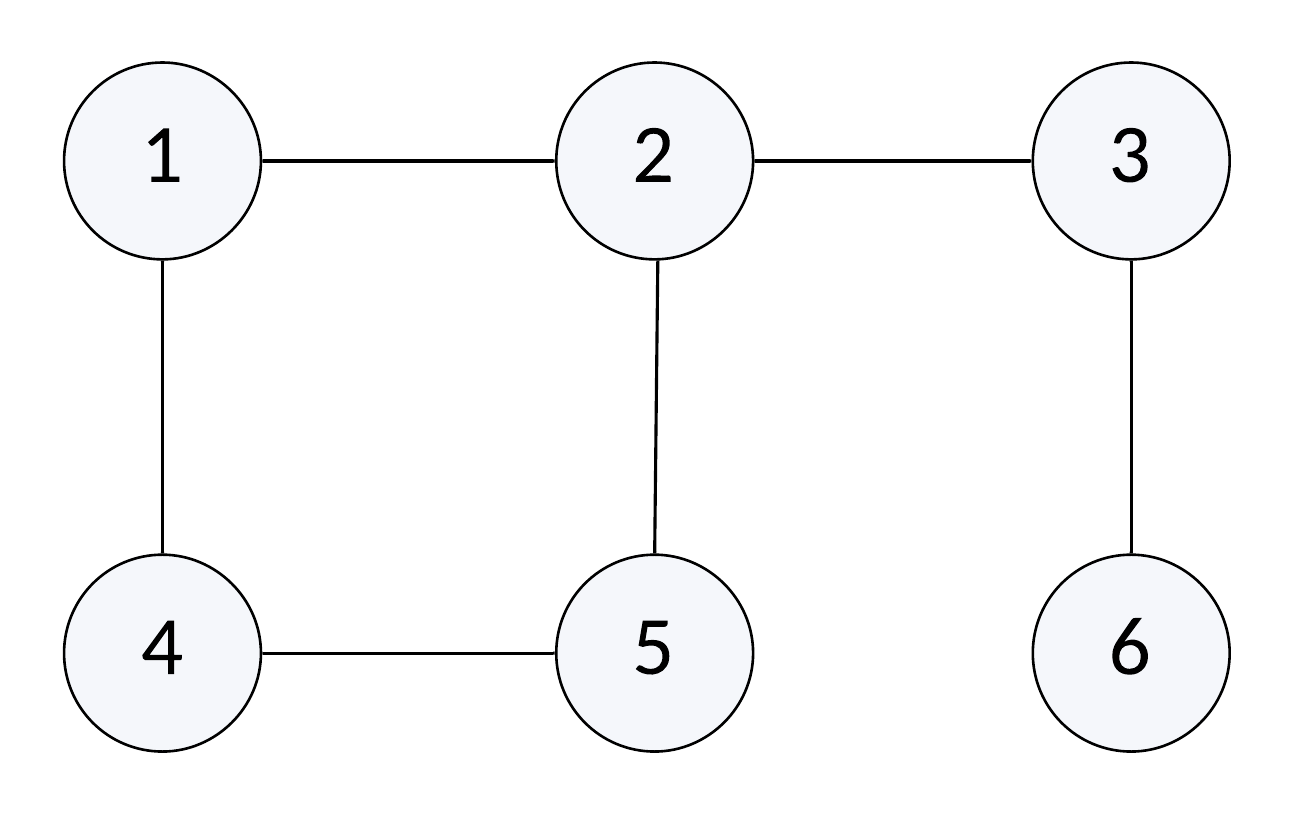} \caption{The communication graph of the 6 agents.}
	\label{ch5a}
\end{figure}

The communication graph of the 6 agents is shown in Fig. \ref{ch5a}.

 Let $\mathbf{x}_i\in\mathbb{R}^6$ be agent $i$'s estimation on the optimal solution. Similar to Section \ref{simu1}, it can be verified that the optimization problem in (\ref{systemsim}) satisfies  Assumptions \ref{a1dis}-\ref{boundednessgra}.
 
 The initial value of $\mathbf{x}$ is selected as $ \mathbf{x}(0)=[2,1,0,3,0,1, 1,1,0,3,0 ,4, 2, 1, 0, 1, 0, 1, 2, 1, 0, 3, 0, 2, 2, 1, 0,$ $ 3, 0, 1, 2, 1, 0, 0, 0, 1]^T$, and  the initial value of $\dot{\mathbf{x}}$ is zero.
 
\begin{figure}[htb]
 	\centering
 	\includegraphics[width=7cm]{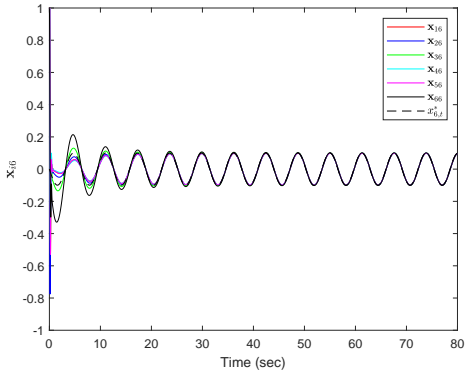} \caption{The agents' estimation on the optimal solution $x_{6,t}^*$.}
 	\label{ch2_dis_allagent}
 \end{figure}

 Fig. \ref{ch2_dis_allagent} shows all agents' estimation on the optimal solution $x_{6,t}^*$. It can be seen that the agents follow the exact optimal solution with high accuracy.
 
\subsubsection{Static Regret}
 
\begin{figure}[htb]
 	\centering
 	\includegraphics[width=7cm]{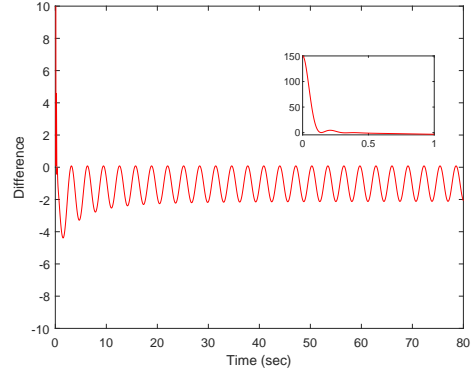} \caption{The difference $\frac{1}{6}\sum_{j=1}^6\sum_{i=1}^6 f_{i,t}(\mathbf{x}_j)- \sum_{i=1}^6f_{i,t}({\tilde{x}}(20))$.}
 	\label{ch22_dis}
 \end{figure}

 Let (\ref{algorithm3}) be the updating law with (\ref{f}), (\ref{eq}), (\ref{fff}), (\ref{pb}) and (\ref{pm}), and the parameters being selected as: $m=-2, \sigma=2$, $b_0=2$, $k_1=2$. Fig. \ref{ch22_dis} shows the difference of the integrated function in the static regret for $T=20$, i.e., $\frac{1}{6}\sum_{j=1}^6\sum_{i=1}^6 f_{i,t}(\mathbf{x}_j)- \sum_{i=1}^6f_{i,t}({\tilde{x}}(20))$, where
 \begin{eqnarray}
 \tilde{x}(20)=\mathop{\text{argmin}}\limits_{x\in\mathbb{R}^6}\int_{0}^{20} f_t(x)dt. \label{simu_offline}
 \end{eqnarray}

 It can be seen that the static regret (the integration of the curve) is upper bounded, which verifies Theorem \ref{theorem4}.
 
\subsubsection{Dynamic Regret}

\begin{figure}[htb]
 	\centering
 	\includegraphics[width=7cm]{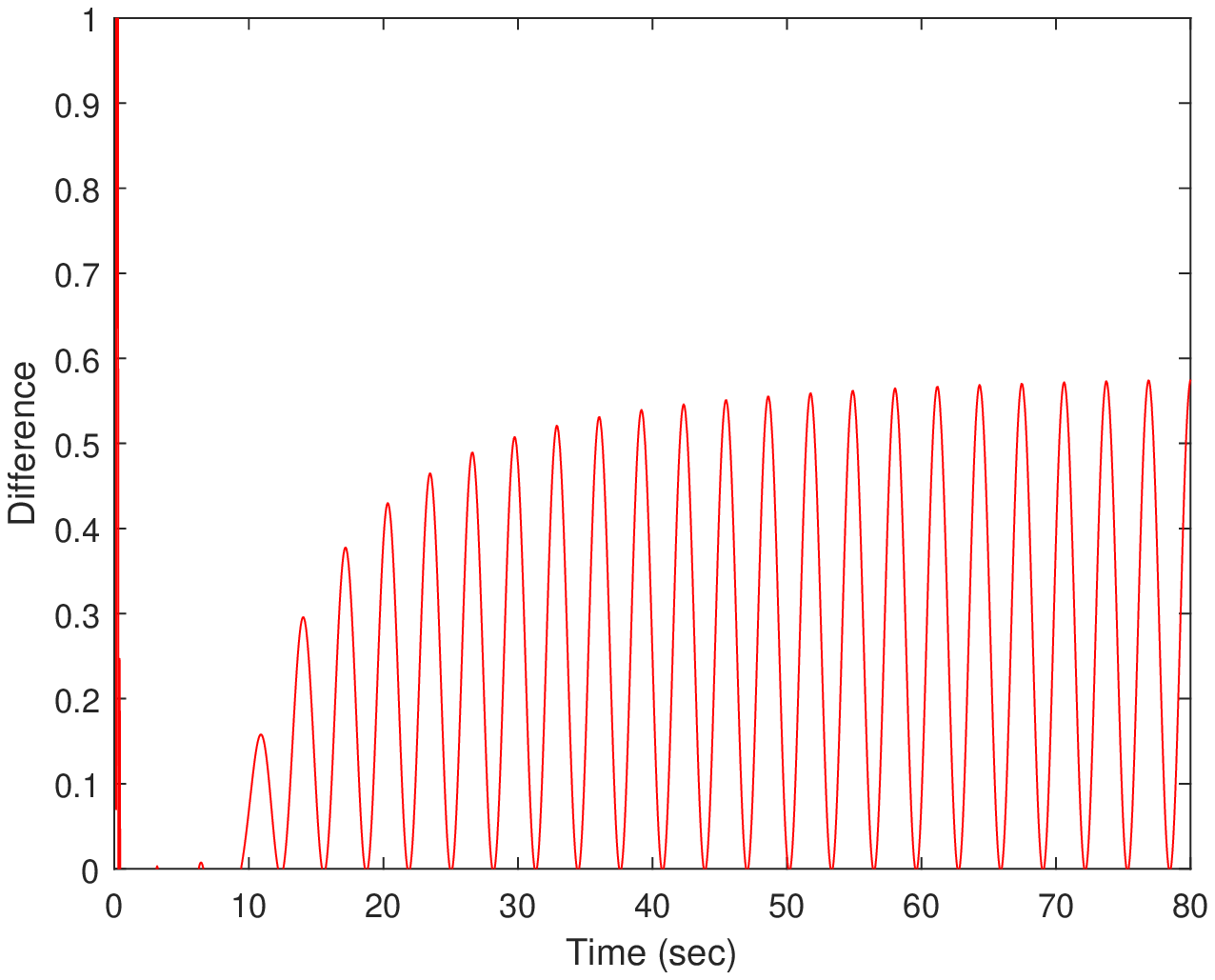} \caption{The difference $\frac{1}{6}\sum_{j=1}^6\sum_{i=1}^6 f_{i,t}(\mathbf{x}_j)- \sum_{i=1}^6f_{i,t}({x}^*_t)$ for $e^{\alpha_t+\beta_t}-\dot{\beta}_te^{\beta_t}=\sigma$.}
 	\label{ch25_dis_normal}
 \end{figure}

\begin{figure}[htb]
	\centering
	\includegraphics[width=7cm]{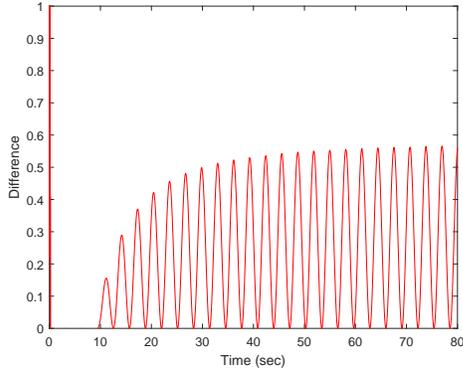} \caption{The difference $\frac{1}{6}\sum_{j=1}^6\sum_{i=1}^6 f_{i,t}(\mathbf{x}_j)- \sum_{i=1}^6f_{i,t}({x}^*_t)$ for $e^{\alpha_t+\beta_t}-\dot{\beta}_te^{\beta_t}=\sigma(t+b_0)$.} \label{ch25_dis_innormal}
	\end{figure}

 Let (\ref{algorithm3}) be the updating law  with (\ref{f}), (\ref{eq}), (\ref{fff}), (\ref{pb}) and (\ref{pm}), and the parameters being selected as: $m=-2, \sigma=2$, $b_0=2$, $k_1=2$. Fig. \ref{ch25_dis_normal} shows the difference of the integrated function in the dynamic regret for $T=80$, i.e., $\frac{1}{6}\sum_{j=1}^6\sum_{i=1}^6 f_{i,t}(\mathbf{x}_j)- \sum_{i=1}^6f_{i,t}({x}^*_t)$, where $e^{\alpha_t+\beta_t}-\dot{\beta}_te^{\beta_t}=\sigma$. Fig. \ref{ch25_dis_innormal} shows the difference when $e^{\alpha_t+\beta_t}-\dot{\beta}_te^{\beta_t}=\sigma(t+b_0)$. It can be seen that the performance is not better when $e^{\alpha_t+\beta_t}-\dot{\beta}_te^{\beta_t}=\sigma(t+b_0)$ as the centralized algorithm in Fig. \ref{ch25}.

\section{Conclusions} \label{S5}

This paper studied the online convex optimization problem by using the proposed OCT-NAG and DOCT-NAG. The online Bregman Lagrangians can generate a family of online optimization algorithms with different scaling conditions. It was shown that for some scaling conditions and under some assumptions, the algorithm achieves a constant static regret and an $O(T)$ dynamic regret. The algorithm was further applied to solve a distributed online optimization problem. Comparable static and dynamic regrets were obtained. In future, we will work on relaxing the assumptions, discretizing the algorithm, extending the work to stochastic settings and solving constrained optimization problems.

\bibliographystyle{IEEEtran}
\bibliography{bib}

\begin{thebibliography}{10}
\providecommand{\url}[1]{#1}
\csname url@samestyle\endcsname
\providecommand{\newblock}{\relax}
\providecommand{\bibinfo}[2]{#2}
\providecommand{\BIBentrySTDinterwordspacing}{\spaceskip=0pt\relax}
\providecommand{\BIBentryALTinterwordstretchfactor}{4}
\providecommand{\BIBentryALTinterwordspacing}{\spaceskip=\fontdimen2\font plus
\BIBentryALTinterwordstretchfactor\fontdimen3\font minus
  \fontdimen4\font\relax}
\providecommand{\BIBforeignlanguage}[2]{{%
\expandafter\ifx\csname l@#1\endcsname\relax
\typeout{** WARNING: IEEEtran.bst: No hyphenation pattern has been}%
\typeout{** loaded for the language `#1'. Using the pattern for}%
\typeout{** the default language instead.}%
\else
\language=\csname l@#1\endcsname
\fi
#2}}
\providecommand{\BIBdecl}{\relax}
\BIBdecl

\bibitem{boyd2004convex}
S.~Boyd and L.~Vandenberghe, \emph{Convex optimization}.\hskip 1em plus 0.5em
  minus 0.4em\relax Cambridge university press, 2004.

\bibitem{bertsekas1997nonlinear}
D.~P. Bertsekas, ``Nonlinear programming,'' \emph{Journal of the Operational
  Research Society}, vol.~48, no.~3, pp. 334--334, 1997.

\bibitem{nesterov27method}
Y.~Nesterov, ``A method of solving a convex programming problem with
  convergence rate $o (\frac{1}{k^2})$,'' in \emph{Sov. Math. Dokl}, vol.~27,
  no.~2.

\bibitem{nesterov2013introductory}
------, \emph{Introductory lectures on convex optimization: A basic
  course}.\hskip 1em plus 0.5em minus 0.4em\relax Springer Science \& Business
  Media, 2013, vol.~87.

\bibitem{qu2017accelerated}
G.~Qu and N.~Li, ``Accelerated distributed nesterov gradient descent,''
  \emph{arXiv preprint arXiv:1705.07176}, 2017.

\bibitem{tatarenko2018accelerated}
T.~Tatarenko, W.~Shi, and A.~Nedi{\'c}, ``Accelerated gradient play algorithm
  for distributed nash equilibrium seeking,'' in \emph{2018 IEEE Conference on
  Decision and Control (CDC)}, pp. 3561--3566.

\bibitem{su2014differential}
W.~Su, S.~Boyd, and E.~Candes, ``A differential equation for modeling
  nesterov’s accelerated gradient method: Theory and insights,'' in
  \emph{Advances in Neural Information Processing Systems}, 2014, pp.
  2510--2518.

\bibitem{pnas}
A.~Wibisono, A.~C. Wilson, and M.~I. Jordan, ``A variational perspective on
  accelerated methods in optimization,'' \emph{proceedings of the National
  Academy of Sciences}, vol. 113, no.~47, pp. E7351--E7358, 2016.

\bibitem{vassilis2018differential}
A.~Vassilis, A.~Jean-Fran{\c{c}}ois, and D.~Charles, ``The differential
  inclusion modeling fista algorithm and optimality of convergence rate in the
  case $b\leq 3$,'' \emph{SIAM Journal on Optimization}, vol.~28, no.~1, pp.
  551--574, 2018.

\bibitem{hazan2016introduction}
E.~Hazan \emph{et~al.}, ``Introduction to online convex optimization,''
  \emph{Foundations and Trends{\textregistered} in Optimization}, vol.~2, no.
  3-4, pp. 157--325, 2016.

\bibitem{popkov2005gradient}
A.~Y. Popkov, ``Gradient methods for nonstationary unconstrained optimization
  problems,'' \emph{Automation and Remote Control}, vol.~66, no.~6, pp.
  883--891, 2005.

\bibitem{Chao}
C.~Sun, M.~Ye, and G.~Hu, ``Distributed time-varying quadratic optimization for
  multiple agents under undirected graphs,'' \emph{IEEE Transactions on
  Automatic Control}, vol.~62, no.~7, pp. 3687--3694, 2017.

\bibitem{rahili2017distributed}
S.~Rahili and W.~Ren, ``Distributed continuous-time convex optimization with
  time-varying cost functions,'' \emph{IEEE Transactions on Automatic Control},
  vol.~62, no.~4, pp. 1590--1605, 2017.

\bibitem{lee2016distributed}
S.~Lee, A.~Ribeiro, and M.~M. Zavlanos, ``Distributed continuous-time online
  optimization using saddle-point methods,'' in \emph{2016 IEEE 55th Conference
  on Decision and Control (CDC)}, pp. 4314--4319.

\bibitem{zinkevich2003online}
M.~Zinkevich, ``Online convex programming and generalized infinitesimal
  gradient ascent,'' in \emph{Proceedings of the 20th International Conference
  on Machine Learning (ICML-03)}, 2003, pp. 928--936.

\bibitem{hazan2007logarithmic}
E.~Hazan, A.~Agarwal, and S.~Kale, ``Logarithmic regret algorithms for online
  convex optimization,'' \emph{Machine Learning}, vol.~69, no. 2-3, pp.
  169--192, 2007.

\bibitem{nedic2010constrained}
A.~Nedic, A.~Ozdaglar, and P.~A. Parrilo, ``Constrained consensus and
  optimization in multi-agent networks,'' \emph{IEEE Transactions on Automatic
  Control}, vol.~55, no.~4, pp. 922--938, 2010.

\bibitem{liu2017convergence}
S.~Liu, Z.~Qiu, and L.~Xie, ``Convergence rate analysis of distributed
  optimization with projected subgradient algorithm,'' \emph{Automatica},
  vol.~83, pp. 162--169, 2017.

\bibitem{shahrampour2017distributed}
S.~Shahrampour and A.~Jadbabaie, ``Distributed online optimization in dynamic
  environments using mirror descent,'' \emph{IEEE Transactions on Automatic
  Control}, vol.~63, no.~3, pp. 714--725, 2017.

\bibitem{hosseini2013online}
S.~Hosseini, A.~Chapman, and M.~Mesbahi, ``Online distributed optimization via
  dual averaging,'' in \emph{52nd IEEE Conference on Decision and
  Control}.\hskip 1em plus 0.5em minus 0.4em\relax IEEE, 2013, pp. 1484--1489.

\bibitem{mateos2014distributed}
D.~Mateos-Nunez and J.~Cort{\'e}s, ``Distributed online convex optimization
  over jointly connected digraphs,'' \emph{IEEE Transactions on Network Science
  and Engineering}, vol.~1, no.~1, pp. 23--37, 2014.

\bibitem{koppel2015saddle}
A.~Koppel, F.~Y. Jakubiec, and A.~Ribeiro, ``A saddle point algorithm for
  networked online convex optimization,'' \emph{IEEE Transactions on Signal
  Processing}, vol.~63, no.~19, pp. 5149--5164, 2015.

\bibitem{hosseini2016online}
S.~Hosseini, A.~Chapman, and M.~Mesbahi, ``Online distributed convex
  optimization on dynamic networks,'' \emph{IEEE Transactions on Automatic
  Control}, vol.~61, no.~11, pp. 3545--3550, 2016.

\bibitem{minimax}
E.~Takimoto and W.~Manfred, ``The minimax strategy for gaussian density
  estimation,'' in \emph{Proc. 13th Annu. Conference on Comput. Learning
  Theory}, 2000, pp. 100--106.

\bibitem{7044563}
E.~C. {Hall} and R.~M. {Willett}, ``Online convex optimization in dynamic
  environments,'' \emph{IEEE Journal of Selected Topics in Signal Processing},
  vol.~9, no.~4, pp. 647--662, June 2015.

\bibitem{besbes}
O.~Besbes, G.~Yonatan, and Z.~Assaf, ``Non-stationary stochastic
  optimization,'' \emph{Operations Research}, vol.~63, no.~5, pp. 1227--1244,
  2015.

\bibitem{NIPS2009_3817}
C.~Hu, W.~Pan, and J.~T. Kwok, ``Accelerated gradient methods for stochastic
  optimization and online learning,'' in \emph{Advances in Neural Information
  Processing Systems 22}, 2009, pp. 781--789.

\bibitem{RenTAC05}
W.~Ren and R.~W. Beard, ``Consensus seeking in multiagent systems under
  dynamically changing interaction topologies,'' \emph{IEEE Transactions on
  Automatic Control}, vol.~50, no.~5, pp. 655--661, 2005.

\bibitem{simonetto2017time}
A.~Simonetto, ``Time-varying convex optimization via time-varying averaged
  operators,'' \emph{arXiv preprint arXiv:1704.07338}, 2017.

\bibitem{duchi2011adaptive}
J.~Duchi, E.~Hazan, and Y.~Singer, ``Adaptive subgradient methods for online
  learning and stochastic optimization,'' \emph{Journal of Machine Learning
  Research}, vol.~12, no. Jul, pp. 2121--2159, 2011.

\bibitem{shalev2012online}
S.~Shalev-Shwartz \emph{et~al.}, ``Online learning and online convex
  optimization,'' \emph{Foundations and Trends{\textregistered} in Machine
  Learning}, vol.~4, no.~2, pp. 107--194, 2012.

\bibitem{hazan2006efficient}
E.~Hazan and S.~Arora, \emph{Efficient algorithms for online convex
  optimization and their applications}.\hskip 1em plus 0.5em minus 0.4em\relax
  Princeton University Princeton, 2006.

\bibitem{Gharesifard14}
B.~Gharesifard and J.~Cort{\'e}s, ``Distributed continuous-time convex
  optimization on weight-balanced digraphs,'' \emph{IEEE Transactions on
  Automatic Control}, vol.~59, no.~3, pp. 781--786, 2014.

\bibitem{liu2014continuous}
S.~Liu, Z.~Qiu, and L.~Xie, ``Continuous-time distributed convex optimization
  with set constraints,'' \emph{IFAC Proceedings Volumes}, vol.~47, no.~3, pp.
  9762--9767, 2014.

\end{thebibliography}

\end{document}